\documentclass[twoside]{mystyle}

\newtheorem{theorem}{Theorem}[subsection]
\newtheorem{corollary}[theorem]{Corollary}
\newtheorem{lemma}[theorem]{Lemma}
\newtheorem{proposition}[theorem]{Proposition}

\newcommand{\qed}{\mbox{$\Box$}}
\renewcommand{\qed}{\hskip0.25em\raisebox{0.6ex}{\framebox[0.5em][l]{\ }}\vspace{.5pc}}
\newcommand{\onespace}{\renewcommand{\baselinestretch}{1}\normalsize}
\newcommand{\twospace}{\renewcommand{\baselinestretch}{1.5}\normalsize}
\newcommand{\threespace}{\renewcommand{\baselinestretch}{1.5}\normalsize}
\newcommand{\num}{\renewcommand{\theequation}{\thetheorem}\addtocounter{theorem}{1}}

\def\theequation{\thetheorem}

\NONUMBIB

\begin{document}

\twospace

\SPECFNSYMBOL{}{}{}{}{}{}{}{}{}%

\AOPMAKETITLE

\AOPAMS{Primary 60B15, 60J15; secondary 20E22.}
\AOPKeywords{Random walk, Markov chain, wreath product, group, Fourier transform, eigenvalue, comparison technique.}
\AOPtitle{RANDOM WALKS ON WREATH PRODUCTS OF GROUPS}
\AOPauthor{Clyde~H.~Schoolfield,~Jr.}
\AOPaffil{Harvard University}
\AOPlrh{C.H. SCHOOLFIELD, JR.}
\AOPrrh{RANDOM WALKS ON WREATH PRODUCTS OF GROUPS}
\AOPAbstract{We bound the rate of convergence to uniformity for certain random walks on the complete monomial groups $G~\wr~S_n$
for any group $G$.  These results provide rates of convergence for random walks on a number of groups of interest: the 
hyperoctahedral group $\mathbb{Z}_2~\wr~S_n$, the generalized symmetric group $\mathbb{Z}_m~\wr~S_n$, and $S_m~\wr~S_n$.
These results provide benchmarks to which many other random walks, modeling a wide range of phenomena, may be compared using
the comparison technique, thereby yielding bounds on the rates of convergence to uniformity for previously intractable random walks.}

\maketitle

\BACKTONORMALFOOTNOTE{3}

\thispagestyle{empty}

\twospace

\section{Introduction.} \label{1}

How many steps does it take for a deck of $n$ cards to achieve near-randomness if, at each step, two randomly chosen cards are
transposed?  This question was answered by Diaconis and Shahshahani (1981).  Now suppose that at each step, the two cards are not 
only transposed, but also possibly flipped.  Suppose that, instead of $n$ cards, there are $n$ wheels and that, at each step, two
wheels are transposed and then spun.  Or suppose that, instead of $n$ wheels, there are $n$ decks of cards and that, at each step, 
two decks are transposed and then shuffled.  How many steps does it take for these processes to achieve near-randomness? 

Or perhaps we are interested in a process that is (at least somewhat) similar in form to one of those described above.  Can we
determine how many steps are needed for it to achieve near-randomness?  These are the types of questions that we consider.

For a certain random walk on the symmetric group $S_n$ that is generated by random transpositions, Diaconis and Shahshahani (1981)
obtained bounds on the rate of convergence to uniformity using group representation theory.  Similarly, we bound the rate of
convergence to uniformity for a random walk on the complete monomial group $G~\wr~S_n$ that is generated by random transpositions,
followed by independent randomizations of the transposed elements.  Specifically, we determine that $\frac{1}{2} n \log n +
\frac{1}{4} n \log (|G|-1|)$ steps are both necessary and sufficient for $\ell^2$ distance to become small.  We also determine that
$\frac{1}{2} n \log n$ steps are both necessary and sufficient for total variation distance to become small.  These results provide
rates of convergence for random walks on a number of groups of interest: the hyperoctahedral group $\mathbb{Z}_2~\wr~S_n$, the
generalized symmetric group $\mathbb{Z}_m~\wr~S_n$, and $S_m~\wr~S_n$.  In the special case of the hyperoctahedral group, the two 
rates of convergence are the same in both metrics.  We also examine a slight variant of this random walk, establishing upper and 
lower bounds on its rate of convergence to uniformity.

The comparison technique was introduced by Diaconis and Saloff-Coste (1993) as a method for bounding the rate of convergence to 
uniformity of a symmetric random walk on a finite group by comparing it to a benchmark random walk whose rate of convergence is 
known.  Our random walks on the hyperoctahedral, generalized symmetric, and complete monomial groups provide such benchmarks to
which many other random walks, modeling a wide range of phenomena, may be compared, thereby yielding bounds on the rates of
convergence to uniformity for previously intractable random walks.  Schoolfield (1998) used the results of this paper to analyze 
two specific examples using this technique, one of which has applications in mathematical biology.  Schoolfield (1998) further 
specialized the comparison technique to random walks generated by random transpositions along the edges of a graph and analyzed 
several examples.

Section~\ref{2} presents the basic properties and results from group theory and representation theory necessary to analyze random
walks on groups.  Section~\ref{3} extends the idea of random transpositions of $n$ cards, analyzed in Section~\ref{2}, to a set of $n$
decks of $m$ cards each and beyond.  

\section{Groups, Representations, and Random Walks.} \label{2}

\subsection{Introduction.} \label{2.1}

Imagine $n$ cards, labeled 1 through $n$, on a table in sequential order.  Independently choose two integers $p$ and $q$ uniformly
from $\{1, 2, \ldots, n\}$.  If $p \neq q$, transpose the cards in positions $p$ and $q$; we denote this by $\tau \in S_n$, where
$\tau$ is the transposition $(p\ q)$.  We shall refer to this procedure as a \emph{random transposition} of $n$ cards.  If $p = q$
(which occurs with probability $1/n$), leave the cards in their current positions.  This action is of course the identity 
permutation, which is denoted by $e \in S_n$. 

If this process is repeated many times, the cards will appear to be in uniformly random order, that is, will appear to be the
result of a \emph{random permutation}.  This process may be modeled formally using a probability measure $P$ (which we shall regard
as a probability mass function) on the symmetric group $S_n$, namely,
\num \begin{equation} \label{2.1.1}
\begin{array}{rcll}

P(e) & := & \displaystyle \frac{1}{n} & \mbox{where $e$ is the identity element}, \vspace{.5pc} \\

P(\tau) & := & \displaystyle \frac{2}{n^2} & \mbox{where $\tau$ is any transposition}, \vspace{.5pc} \\

P(\pi) & := & 0 & \mbox{otherwise}.

\end{array}
\end{equation}

\noindent
Repeating the process above $k$ times is modeled as the convolution $P^{*k}$ of the measure $P$ with itself $k$ times.  Since there
are $n!$ elements in $S_n$, the uniform probability measure on the set of all permutations of $S_n$ is given by
\num \begin{equation} \label{2.1.2}
U(\pi) \ \ := \ \ \frac{1}{n!} \ \ \ \mbox{for every $\pi \in S_n$}.
\end{equation}

The following result, which is Theorem 1 of Diaconis and Shahshahani (1981) and was later included in Diaconis (1988) as Theorem 5
in Section D of Chapter 3, establishes an upper bound on both the total variation distance and the $\ell^2$ distance (both defined 
in Section~\ref{2.6}) between $P^{*k}$ and $U$.

\begin{theorem} \label{2.1.3}
Let $P$ and $U$ be the probability measures on the symmetric group $S_n$ defined in (\ref{2.1.1}) and (\ref{2.1.2}),
respectively.  Let $k = \frac{1}{2} n \log n + cn$.  Then there exists a universal constant $a > 0$ such that

\[ \| P^{*k} - U \|_{\mbox{\rm \scriptsize TV}} \ \ \leq \ \ \mbox{$\frac{1}{2}$} \left(n!\right)^{1/2} \| P^{*k} - U \|_2 \ \
\leq \ \ ae^{-2c} \ \ \ \mathrm{for\ all\ } \mbox{$c > 0$}. \]
  
\end{theorem}

In the following sections we present the results that were needed to prove this theorem and which are used to prove analogous 
results in later sections.  In Section~\ref{2.2} we present basic definitions from group theory.  In Section~\ref{2.3} we present
basic results from the theory of group representations, while Section~\ref{2.4} concentrates on the characters of these 
representations.  In Section~\ref{2.5} we introduce the Fourier transform, and we show how it may be used to bound the distance to
uniformity of a random walk on a group in Section~\ref{2.6}.  In Section~\ref{2.7} we show how the results from these previous
sections were applied to the random walk on the symmetric group defined above, proving Theorem~\ref{2.1.3} and a matching lower bound.

\subsection{Group Theory.} \label{2.2}

We now present basic definitions from the theory of groups which will be needed to conduct our analysis.  A more detailed 
introduction to this subject may be found in Chapter 1 of Alperin and Bell (1995).

A \emph{group} is a non-empty set $G$ with a binary operation on $G$, usually called multiplication, which satisfies the following
three axioms:  (i) Multiplication is associative, (ii) There is a unique \emph{identity} element $e \in G$, and (iii) For every $g
\in G$ there is a unique \emph{inverse} element $g^{-1} \in G$ such that $gg^{-1} = e = g^{-1}g$.  The number of elements of a
finite group $G$ is called the \emph{order} of $G$ and is denoted by $|G|$.  

A subset $H$ of a group $G$ is called a \emph{subgroup} of $G$ if it forms a group under the multiplication of $G$ restricted to
$H$.  A subgroup $N$ of $G$ is called a \emph{normal subgroup} of $G$ if $gNg^{-1} \subseteq N$ for all $g \in G$, where $gNg^{-1} 
:= \{ gng^{-1} \in G  : n \in N \}$.  If $gh = hg$ for all $g,h \in G$, then $G$ is called an \emph{abelian} group.  Notice that
every subgroup of an abelian group is a normal subgroup.

The set $kH := \{ kh \in G  : h \in H \}$, where $k \in G$ and $H$ is a subgroup of $G$, is called a \emph{left coset} of $H$ in
$G$.  (Right cosets are defined analogously.)  The set of all left cosets of $H$ in $G$ is called the \emph{left coset space} and
is denoted by $G/H$.  If $|G/H| = n$, a set of \emph{left coset representatives} $\{ k_1, k_2, \ldots, k_n \}$ may be chosen so
that the left cosets $k_1H, k_2H, \ldots, k_nH$ comprise precisely the coset space $G/H$.  If $N$ is a normal subgroup of $G$, then
$G/N$ is also a group known as the \emph{quotient group}.

One group of interest to us is the \emph{cyclic group} $\mathbb{Z}_n$, which is the set $\{0, 1, 2, \ldots, n-1\}$ under the operation
of addition mod $n$.  Notice that $\mathbb{Z}_n$ is abelian and that $|\mathbb{Z}_n| = n$.  Another group of interest to us is the
\emph{symmetric group} $S_n$, which is the set of all permutations of $\{1, 2, \ldots, n\}$ under the operation of composition of
functions (where we compose functions from right to left).  Notice that $S_n$ is \emph{not} abelian for any $n \geq 3$, and that
$|S_n| = n!$.

Two elements $g$ and $h$ of $G$ are called \emph{conjugate} if there exists some $k \in G$ such that $h = kgk^{-1}$.  The set of
elements conjugate to a particular element $g \in G$ form the \emph{conjugacy class} of $g$.  Conjugacy is an equivalence
relation and partitions a group $G$ into disjoint subsets, namely, the conjugacy classes.

The set of ordered pairs of elements of groups $G_1$ and $G_2$ with multiplication defined componentwise is called the \emph{direct
product} of $G_1$ and $G_2$ and is denoted by $G_1~\times~G_2$.  The group of elements $(g_1, g_2, \ldots, g_n; \pi) \in G^n \times
S_n$, with multiplication defined by 

\[ (h_1, \ldots, h_n; \sigma) \cdot (g_1, \ldots, g_n; \pi) = (h_1 \cdot g_{\sigma^{-1}(1)}, \ldots, h_n \cdot g_{\sigma^{-1}(n)};
\sigma \pi), \] 

\noindent
is called the \emph{wreath product} of $G$ with $S_n$ and is denoted $G~\wr~S_n$.  The identity element of $G~\wr~S_n$ is $(e,
\ldots, e; e)$, and in $G~\wr~S_n$ we have $(g_1, \ldots, g_n; \pi)^{-1} = (g_{\pi(1)}^{-1}, \ldots, g_{\pi(n)}^{-1}; \pi^{-1})$.
This definition is equivalent to saying that $G~\wr~S_n$ is the \emph{semidirect product} of $G^n$ with $S_n$, where the action of
$S_n$ on $G^n$ is $\pi(g_1, g_2, \ldots, g_n) = (g_{\pi^{-1}(1)}, g_{\pi^{-1}(2)}, \ldots, g_{\pi^{-1}(n)})$.  For a general
definition of semidirect product see Alperin and Bell (1995) or Simon (1996).  Three special cases of particular interest to us are
the \emph{hyperoctahedral group} $\mathbb{Z}_2~\wr~S_n$, the \emph{generalized symmetric group} $\mathbb{Z}_m~\wr~S_n$, and $S_m \wr
S_n$.  

\subsection{Representation Theory.} \label{2.3}

We now present basic properties and results from the representation theory of groups which will be needed to conduct our
analysis.  A more detailed introduction to this subject may be found in Chapters 1 through 3 of Serre (1977) or in Chapter 2 of
Diaconis (1988).  Other sources include Alperin and Bell (1995) and Simon (1996).

Suppose that $V$ is a finite-dimensional vector space over the complex numbers.  The \emph{general linear group} \textbf{GL}$(V)$ is
the group of isomorphisms of $V$ onto itself.  An element of \textbf{GL}$(V)$ is a linear mapping of $V$ onto $V$.  Such a map has
an inverse which is also linear. 

A \emph{representation} of a finite group $G$ in $V$ is a homomorphism $\rho: G \longrightarrow\ $\textbf{GL}$(V)$.  The choice of a
basis for $V$ assigns an invertible matrix $\rho(g)$ to each $g \in G$.  There exists $\rho(g) \in\ $\textbf{GL}$(V)$ for each $g
\in G$ such that $\rho(g_1 g_2) = \rho(g_1) \rho(g_2)$ for all $g_1, g_2 \in G$.  This implies that $\rho(e) = I$ and that
$\rho(g^{-1}) = \rho(g)^{-1}$. 

Two representations of $G$, say $\rho_1$ in $V_1$ and $\rho_2$ in $V_2$, are said to be \emph{isomorphic} (or \emph{equivalent}) if 
there exists a linear isomorphism $\mu: V_1 \longrightarrow V_2$ such that $\mu \circ \rho_1(g) = \rho_2(g) \circ \mu$ for all $g
\in G$.  The \emph{dimension} (or \emph{degree}) of $\rho$ is defined to be the dimension of $V$ and is denoted by $d_{\rho}$.  
Since we shall be interested in representations only up to equivalence, we may without loss of generality assume that $V =  
\mathbb{C}^n$ when $d_{\rho} = n$.  The \emph{trivial representation} is the one-dimensional representation with $V = \mathbb{C}$
that sends every element of $G$ to 1. 

A vector subspace $W \subseteq V$ is \emph{stable} under $G$ if $\rho(g)w \in W$ for all $w \in W$ and all $g \in G$.  If $W
\subseteq V$ is stable under $G$, then $\rho$ restricted to \textbf{GL}$(W)$ gives a \emph{subrepresentation} of $\rho$.  A 
representation $\rho$ is \emph{irreducible} if the only vector subspaces of $V$ which are stable under $G$ are $V$ itself and the
trivial subspace.  An important relationship between irreducible representations and conjugacy classes is given by the following,
which is Theorem 7 in Section 2.5 of Serre (1977).

\begin{proposition} \label{2.3.1}
The number of nonisomorphic irreducible representations of $G$ equals the number of conjugacy classes of $G$.
\end{proposition}

An important relationship between the dimensions of the irreducible representations of $G$ and the order of $G$ is given by the
following, which is Corollary 2(a) of Proposition 5 in Section 2.4 of Serre (1977).

\begin{lemma} \label{2.3.2}
Suppose that $\{ \rho_1, \rho_2, \ldots, \rho_s \}$ is a complete set of nonisomorphic irreducible representations of $G$ with
dimensions $d_{\rho_1}, d_{\rho_2}, \ldots, d_{\rho_s}$, respectively.  Then

\[ \sum_{i=1}^s d_{\rho_i}^2 \ \ = \ \ |G|. \]

\end{lemma}

It follows directly from Proposition~\ref{2.3.1} and Lemma~\ref{2.3.2} that a group $G$ is abelian if and only if all of its
irreducible representations are one-dimensional.

The following useful result is \emph{Schur's Lemma}, which is Proposition 4 in Section 2.2 of Serre (1977).

\begin{lemma} \label{2.3.3}
Suppose that $\rho_1: G \longrightarrow\ $ \emph{\textbf{GL}}$(V_1)$ and $\rho_2: G \longrightarrow\ $ \emph{\textbf{GL}}$(V_2)$ are
irreducible representations of $G$ and that $\mu: V_1 \longrightarrow V_2$ is such that $\mu \circ \rho_1(g) = \rho_2(g) \circ
\mu$ for all $g \in G$.  Then \emph{(i)} if $\rho_1$ and $\rho_2$ are not isomorphic, it follows that $\mu = 0;$ and \emph{(ii)} if
$V_1 = V_2$ and $\rho_1 = \rho_2$, it follows that $\mu$ is a constant times the identity.
\end{lemma}

The \emph{direct sum} $A \oplus B$ of an $m~\times~m$ matrix $A$ and an $n~\times~n$ matrix $B$ is the $(m + n)~\times~(m + n)$
block diagonal matrix \onespace $\left[ \begin{array}{cc} A & 0 \\ 0 & B \\ \end{array} \right]$ \twospace.  The direct sum $\rho_1
\oplus \rho_2$ of two representations is then defined by $\left(\rho_1 \oplus \rho_2\right)(g) = \rho_1(g) \oplus \rho_2(g)$.  By use
of the direct sum, the irreducible representations of $G$ can be used to construct all other representations of $G$, as described in
the following, which is Theorem~2 in Section 1.4 of Serre (1977).

\begin{proposition} \label{2.3.4}
Every representation of $G$ is the direct sum of irreducible representations of $G$.
\end{proposition}

The \emph{tensor product} $A \otimes B$ of an $m~\times~m$ matrix $A$ and an $n~\times~n$ matrix $B$ is an $mn~\times~mn$ matrix
which is constructed in the following manner.  Begin with an $m~\times~m$ block matrix in which each of the $m^2$ blocks is the
matrix $B$.  Then multiply the block in position $(i,j)$ by the scalar $a_{ij} \in A$ for $1 \leq i,j \leq m$.  The tensor product
$\rho_1 \otimes \rho_2$ of two representations is then defined by $\left(\rho_1 \otimes \rho_2\right)(g) = \rho_1(g) \otimes
\rho_2(g)$.  By use of the tensor product, the irreducible representations of $G_1$ and $G_2$ can be used to construct
all the irreducible representations of their direct product, as described in the following, which is Theorem 10 in Section 3.2 of
Serre (1977).

\begin{proposition} \label{2.3.5}
Suppose that $\rho_1$ and $\rho_2$ are irreducible representations of $G_1$ and $G_2$, respectively.  Then $\rho_1 \otimes \rho_2$
is an irreducible representation of $G_1~\times~G_2$.  Furthermore, each irreducible representation of $G_1~\times~G_2$ is
isomorphic to such a representation.
\end{proposition}

Suppose that $H$ is a subgroup of $G$ and that $\rho$ is a representation of $G$.  A representation $\rho \downarrow_H^G$ of $H$,
known as the \emph{restricted representation}, can be constructed from $\rho$ by defining

\[ \rho \downarrow_H^G (h) \ \ := \ \ \rho(h) \ \ \ \mathrm{for\ each\ } \mbox{$h \in H$}. \]

Now suppose that $H$ is a subgroup of $G$ and that $\rho$ is a representation of $H$.  Also suppose that $\{k_1, k_2, \ldots, k_n\}$
is a complete set of left coset representatives of $H$ in $G$.  A representation $\rho \uparrow_H^G$ of $G$, known as the
\emph{induced representation}, can be constructed from $\rho$ by defining, for each $g \in G$, an $n~\times~n$ block matrix $\rho
\uparrow_H^G (g)$ whose block in position $(i,j)$, for $1 \leq i,j \leq n$, is the $d_{\rho}~\times~d_{\rho}$ matrix

\[ \rho \uparrow_H^G (g)_{ij} \ \ := \ \ \left\{ \begin{array}{ll} 
\rho(k_i^{-1} g k_j)  &  \mathrm{if\ } \mbox{$k_i^{-1} g k_j \in H$},  \\
0  &  \mathrm{if\ } \mbox{$k_i^{-1} g k_j \not\in H$}.  \\
\end{array} \right. \]

\noindent
The induced representation does not depend on the choice of coset representatives.

\subsection{Character Theory.} \label{2.4}

The \emph{character} of the representation $\rho$ at the element $g \in G$ is defined to be the trace of $\rho(g)$ and is denoted by
$\chi_{\rho}(g)$.  Notice that the character of a representation is independent of the choice of basis of $V$.  The characters of
the irreducible representations of $G$ are called \emph{irreducible characters}.  The choice of the term ``character'' is to
emphasize that it characterizes the representation; according to Corollary 2 of Theorem 4 in Section 2.3 of Serre (1977), two
representations are isomorphic if and only if they have the same character.

Some important properties of characters are given in the following, which is Proposition 1 in Section 2.1 of Serre (1977).

\begin{lemma} \label{2.4.1}
Suppose that $\chi$ is the character of a representation $\rho$ of $G$ with dimension $d_{\rho}$.  Then
\[ \begin{array}{rll}
\mathrm{(a)}  &  \chi(e) = d_{\rho}  &  \mathrm{for\ } \mbox{$e \in G$},  \\
\mathrm{(b)}  &  \chi(g^{-1}) = \overline{\chi(g)}  &  \mathrm{for\ } \mbox{$g \in G$},  \\
\mathrm{(c)}  &  \chi(hgh^{-1}) = \chi(g)  &  \mathrm{for\ } \mbox{$g,h \in G$}.  \\
\end{array} \]
\end{lemma}

The characters of direct sums and tensor products may be easily calculated by use of the following, which is Proposition 2 in
Section 2.1 of Serre (1977).

\begin{lemma} \label{2.4.2}
Suppose that $\rho_1$ and $\rho_2$ are representations of $G$ with characters $\chi_1$ and $\chi_2$, respectively.  Then the
character of the direct sum $\rho_1 \oplus \rho_2$ is $\chi_1 + \chi_2$ and the character of the tensor product $\rho_1 \otimes
\rho_2$ is $\chi_1 \cdot \chi_2$.
\end{lemma}

The characters of induced representations may be calculated by use of the following, which is Corollary 6 in Section 16 of Alperin
and Bell (1995).

\begin{lemma} \label{2.4.3}
Suppose that $\chi_H$ is the character of an irreducible representation $\rho_H$ of a subgroup $H$ of a finite group $G$.  For $g
\in G$, suppose that the number $t$ of conjugacy classes of $H$ whose members are conjugate in $G$ to $g$ is positive.  Let $h_1,
h_2, \ldots, h_t$ be representatives of these $t$ conjugacy classes of $H$ and let $k_1, k_2, \ldots, k_t$ be the sizes of these
classes.  Let $\ell$ be the size of the conjugacy class of $g$ in $G$.  Then the value at $g$ of the character $\chi$ of the induced
representation $\rho_H \uparrow_H^G$ of $G$ is given by

\[ \chi(g) \ \ = \ \ \frac{|G|}{|H|} \sum_{i=1}^t \frac{k_i}{\ell} \chi_{H}(h_i). \]

\end{lemma}
\vspace{.5pc}

Suppose that $f_1$ and $f_2$ are functions on a group $G$.  The \emph{inner product} of $f_1$ and $f_2$ is defined by

\[ \left\langle f_1 , f_2 \right\rangle_G \ \ := \ \ \frac{1}{|G|} \sum_{g \in G} f_1(g) \overline{f_2(g)}. \]

A function on $G$ which is constant on each conjugacy class of $G$ is called a \emph{class function}.  Notice that 
Lemma~\ref{2.4.1}(c) asserts that the character $\chi$ of a representation $\rho$ of $G$ is a class function.  An important
relationship between the irreducible characters of $G$ and the space of class functions on $G$ is given by the following, which is
Theorem 6 in Section 2.5 of Serre (1977).

\begin{proposition} \label{2.4.4}
The characters $\chi_1, \chi_2, \ldots, \chi_s$ of the irreducible representations of $G$ form an orthonormal basis for the Hilbert
space of class functions on $G$ with respect to the inner product $\left\langle \cdot , \cdot \right\rangle$ defined above.
\end{proposition}

\subsection{The Fourier Transform.} \label{2.5}

We now introduce an extremely important tool for our analysis.  Suppose that $P$ is a function on $G$.  The \emph{Fourier transform}
$\widehat{P}$ of $P$ at the representation $\rho$ is defined to be the matrix

\[ \widehat{P}(\rho) \ \ := \ \ \sum_{g \in G} P(g) \rho(g). \]

Suppose that $P$ and $Q$ are functions on $G$.  The \emph{convolution} $P * Q$ of $P$ and $Q$ is defined, for all $g \in G$, by

\[ P * Q(g) \ \ := \ \ \sum_{h \in G} P(gh^{-1}) Q(h). \]

The Fourier transform converts the convolution of functions $P$ and $Q$ into the (argumentwise) multiplication of their
transforms $\widehat{P}$ and $\widehat{Q}$, i.e.,

\[ \widehat{P * Q}(\rho) \ \ = \ \ \widehat{P}(\rho) \widehat{Q}(\rho). \]

In the special case when $P$ is a class function, it is a consequence of Schur's Lemma (\ref{2.3.3}) that the Fourier transform may
be calculated easily by use of the following, which is Lemma 5 of Diaconis and Shahshahani (1981).

\begin{lemma} \label{2.5.1}
Suppose that $\rho$ is an irreducible representation of a finite group $G$ with character $\chi$ and that $P$ is a class function.
For each conjugacy class $i$, let $P_i$ be the constant value of $P$ on the class, let $n_i$ be the cardinality of the class, and
let $\chi_i$ be the constant value of $\chi$ on the class.  Then the Fourier transform of $P$ is given by

\[ \widehat{P}(\rho) \ \ = \ \ \left[ \frac{1}{d_{\rho}} \sum_{i=1}^s P_i n_i \chi_i \right] I, \]

\noindent
where $d_\rho$ is the dimension of $\rho$, $I$ is the $d_\rho$-dimensional identity matrix, and the sum is taken over distinct
conjugacy classes.  
\end{lemma}

The Fourier transforms of any distribution at the trivial representation and of the uniform distribution at any nontrivial
representation have special forms, as described in the following, which is an immediate consequence of the preceding lemma and
Proposition~\ref{2.4.4}.

\begin{corollary} \label{2.5.2}
Suppose that $P$ is any probability measure and that $U$ is the uniform probability measure, both defined on a finite group $G$.
Then $\widehat{P}(\rho_0) = 1$ at the trivial representation $\rho_0$ of $G$ and $\widehat{U}(\rho)$ is the $d_{\rho} \times
d_{\rho}$ zero matrix for each nontrivial representation $\rho$ of $G$.
\end{corollary}

The function $P$ on $G$ may be reconstructed from its Fourier transform by use of the following, which is the \emph{Fourier
inversion formula} in Section 6.2 of Serre (1977).

\begin{lemma} \label{2.5.3}
Suppose that $P$ is a function on $G$ and that $\widehat{P}$ is its Fourier transform.  Then for each $g \in G$

\[ P(g) \ \ = \ \ \frac{1}{|G|} \sum_{i=1}^s d_{\rho_i} \mbox{\rm tr} \left( \rho_i(g^{-1})\widehat{P}(\rho_i) \right), \]

\noindent
where $\rho_1, \rho_2, \ldots, \rho_s$ are the nonisomorphic irreducible representations of $G$ with dimensions
$d_{\rho_1}, d_{\rho_2}, \ldots, d_{\rho_s}$, respectively.
\end{lemma}

A consequence of this formula is another useful result, which is the \emph{Plancherel formula} in Section 6.2 of Serre (1977).

\begin{lemma} \label{2.5.4}
Suppose that $P$ and $Q$ are functions on $G$ and that $\widehat{P}$ and $\widehat{Q}$ are their Fourier transforms.  Then

\[ \sum_{g \in G} P(g)Q(g^{-1}) \ \ = \ \ \frac{1}{|G|} \sum_{i=1}^s d_{\rho_i} \mbox{\rm tr} \left(
\widehat{P}(\rho_i)\widehat{Q}(\rho_i) \right), \]

\noindent
where $\rho_1, \rho_2, \ldots, \rho_s$ are the nonisomorphic irreducible representations of $G$ with dimensions
$d_{\rho_1}, d_{\rho_2}, \ldots, d_{\rho_s}$, respectively.
\end{lemma}

\subsection{Random Walks on Groups.} \label{2.6}

We now turn our attention to the subject of random walks on groups.  A more detailed introduction to this subject may be found
in Chapter 3 of Diaconis (1988). 

Suppose that $P$ is a probability measure defined on a group $G$.  Let $\xi_1, \xi_2, \ldots$ be a sequence of independent
$G$-valued random variables each distributed according to $P$.  A \emph{random walk} on $G$ is a sequence \textbf{X} $= (X_0, X_1,
X_2, \ldots)$ defined by $X_0 := e \in G$ and $X_n := \xi_n \xi_{n-1} \cdots \xi_1$ for all $n \geq 1$.  Notice that \textbf{X} is a
Markov chain with state space~$G$:

\[ \begin{array}{l}

\mathbb{P}\{ X_n = g_n \ | \ X_0 = e, \ X_1 = g_1, \ \ldots, \ X_{n-1} = g_{n-1} \} \vspace{1pc} \\

\ \ \ = \ \mathbb{P}\{ \xi_n \xi_{n-1} \cdots \xi_1 = g_n \ | \ \xi_1 = g_1, \ \ldots, \ \xi_{n-1} \cdots \xi_1 = g_{n-1} \} 
\vspace{1pc} \\

\ \ \ = \ \mathbb{P}\{ \xi_n = g_n g_{n-1}^{-1} \} \ \ = \ \ P(g_n g_{n-1}^{-1}) 

\end{array} \]

\noindent
for all $n \geq 1$ and all $g_1, \ldots, g_n \in G$.  In this way a probability measure $P$ on $G$ induces a transition matrix
\textbf{P}, where the entry of \textbf{P} at the intersection of the row corresponding to $g \in G$ and the column corresponding to
$h \in G$ is given by

\[ \textbf{P}_{g h} \ \ = \ \ P(hg^{-1}). \]

In the special case where $P$ is a class function, we may determine all the eigenvalues of the transition matrix \textbf{P},
together with their multiplicities, by use of the following, which is Corollary 3 of Diaconis and Shahshahani (1981).

\begin{lemma} \label{2.6.1}
Suppose $P$ is a probability measure defined on a finite group $G$ and that $P$ is a class function.  Let \emph{\textbf{P}} be the
transition matrix of the Markov chain induced by the probability measure $P$.  Then, for each irreducible representation $\rho$ of
$G$, there is an eigenvalue $\pi_{\rho}$ of \emph{\textbf{P}} occurring with algebraic multiplicity $d_{\rho}^2$ such that

\[ \pi_{\rho} \ \ = \ \ \frac{1}{d_{\rho}} \sum_{i=1}^s P_i n_i \chi_i, \]

\noindent
where the sum is taken over distinct conjugacy classes.
\end{lemma}

\noindent
Notice that $\pi_{\rho} \cdot I$ in the lemma above is exactly the value of $\widehat{P}(\rho)$ in Lemma~\ref{2.5.1}.

In order to discuss the convergence of these random walks to their stationary distributions, we need a metric between probability
measures.  Suppose that $P$ and $Q$ are two probability measures defined on a finite group $G$.  The \emph{total variation
distance} between $P$ and $Q$ is defined by

\[ \| P - Q \|_{\mbox{\rm \scriptsize TV}} \ \ := \ \ \max_{A \subseteq G} |P(A) - Q(A)|, \]

\noindent
while the \emph{$\ell^1$ distance} between $P$ and $Q$ is defined by

\[ \| P - Q \|_1 \ \ := \ \ \sum_{g \in G} |P(g) - Q(g)|. \]

\noindent
Notice that $\| P - Q \|_{\mbox{\rm \scriptsize TV}} \ = \ \frac{1}{2} \| P - Q \|_1$.  The \emph{$\ell^2$ distance} between
$P$ and $Q$ is defined by

\[ \| P - Q \|_2 \ \ := \ \ \left( \sum_{g \in G} |P(g) - Q(g)|^2 \right)^{1/2}. \]

\noindent
All three of these measures of distance are indeed metrics.  It is a direct consequence of the Cauchy-Schwarz inequality that

\[ \| P - Q \|_{\mbox{\rm \scriptsize TV}}^2 \ \ \leq \ \ \displaystyle \mbox{$\frac{1}{4}$} |G| \| P - Q \|_2^2. \]

We are now able to bound the distance to uniformity of a probability measure $P$ in terms of its Fourier transform by use of the
following, which is the Upper Bound Lemma in Section B of Chapter 3 of Diaconis (1988).

\begin{lemma} \label{2.6.2}
Suppose that $P$ is a probability measure defined on a finite group $G$.  Then 

\[ \| P - U \|_{\mbox{\rm \scriptsize TV}}^2 \ \ \leq \ \ \mbox{$\frac{1}{4}$} |G| \cdot \| P - U \|_2^2 \ \ = \ \
\mbox{$\frac{1}{4}$} \sum_{\rho} d_{\rho}\ \mbox{\rm tr} \left(\widehat{P}(\rho)\widehat{P}(\rho)^*\right) \]

\noindent
where the sum is taken over all \emph{nontrivial} irreducible representations of $G$ and $\widehat{P}(\rho)^*$ is the conjugate
transpose of $\widehat{P}(\rho)$.
\end{lemma}

\noindent
In the Upper Bound Lemma (\ref{2.6.2}), the inequality follows from the Cauchy-Schwarz inequality; the equality follows from the
Plancherel formula (\ref{2.5.4}) and Corollary~\ref{2.5.2}.

\subsection{Random Walk on the Symmetric Group.} \label{2.7}

We now return to the random walk on the symmetric group defined in Section~\ref{2.1} and the proof of Theorem~\ref{2.1.3}.  We include
only those portions of the proof that will be needed to conduct our analysis in later sections.  A detailed introduction to the
symmetric group and its representation theory may be found in Chapters 1 and 2 of James and Kerber (1981).  Another source is Sagan
(1991).

For any $\pi \in S_n$, there is an associated $n$-dimensional vector $a = (a_1, a_2, \ldots, a_n)$, called the \emph{cycle type} of
$\pi$, where $a_i$ is the numbers of cycle factors of $\pi$ of length $i$ for $1 \leq i \leq n$.  According to Lemma 1.2.6 of James
and Kerber (1981), two elements of $S_n$ are conjugate if and only if their cycle types are identical; hence, there is a one-to-one
correspondence between these cycle type vectors and conjugacy classes of $S_n$.  Thus, since the transpositions in $S_n$ form their
own conjugacy class, the probability measure $P$ defined in (\ref{2.1.1}) is a class function.

A vector $[\lambda] = [\lambda_1, \lambda_2, \ldots, \lambda_k]$ such that  $\lambda_1 \geq \lambda_2 \geq \cdots \geq \lambda_k >
0$ and $\lambda_1 + \lambda_2 + \cdots + \lambda_k = n$ is called a \emph{partition} of $n$.  According to Theorem 2.1.11 of James 
and Kerber (1981), there is a one-to-one correspondence between nonisomorphic irreducible representations of $S_n$ and partitions
$[\lambda]$ of $n$.  Thus we have identified all of the irreducible representations of $S_n$, over which the summation is taken in 
the Upper Bound Lemma (\ref{2.6.2}). 

Lemmas~\ref{2.5.1} and \ref{2.4.1}(a) are then used to calculate the Fourier transform of an irreducible representation
$\rho_{[\lambda]}$ of $S_n$, which is

\[ \widehat{P}(\lambda) = \left[ \frac{1}{n} + \frac{n-1}{n} r(\lambda) \right] I, \]

\noindent
where $r(\lambda) := \chi_{[\lambda]}(\tau)/d_{[\lambda]}$, $\chi_{[\lambda]}(\tau)$ is the character of $\rho_{[\lambda]}$ at any
transposition $\tau$, $d_{[\lambda]}$ is the dimension of $\rho_{[\lambda]}$, and $I$ is the $d_{[\lambda]}$-dimensional identity
matrix.

It then follows from the Upper Bound Lemma (\ref{2.6.2}) that

\[ \| P^{*k} - U \|_{\mbox{\rm \scriptsize TV}}^2 \ \ \leq \ \ \mbox{$\frac{1}{4}$} n! \| P^{*k} - U \|_2^2 \ \ = \ \
\mbox{$\frac{1}{4}$} \sum_{[\lambda]} d_{[\lambda]}^2\ \left[ \frac{1}{n} + \frac{n-1}{n} r(\lambda) \right]^{2k} \]

\noindent
where the term not to be included in the summation occurs when $[\lambda] = [n]$, which corresponds to the trivial representation
of $S_n$.

The following formulas, found is Section D of Chapter 3 and Section B of Chapter~7, respectively, of Diaconis (1988) are used to
calculate the numerical value of the Fourier transform.

\begin{lemma} \label{2.7.1}
Suppose that $\rho$ is an irreducible representation of $S_n$ corresponding to the partition $[\lambda] = [\lambda_1, \ldots,
\lambda_k]$ of $n$.  Let $r(\lambda) := \chi_{[\lambda]}(\tau) / d_{[\lambda]}$ with $\tau \in S_n$.  Then

\[ \begin{array}{rcl}

r(\lambda)  &  =  &  \displaystyle \frac{1}{n(n-1)} \sum_{j=1}^k \left[ \lambda_j^2 \ - \ (2j-1) \lambda_j \right] \ \ \
\mathrm{and} \vspace{1pc} \\

d_{[\lambda]}  &  =  &  \displaystyle n! \ \mbox{\rm det} \displaystyle \left( \frac{1}{(\lambda_i - i + j)!} \right)_{1 \leq 
i,j \leq k},

\end{array} \]

\noindent
with $1/m! := 0$ if $m < 0$.

\end{lemma}

A lengthy, detailed discussion in Section D of Chapter 3 of Diaconis (1988) determines the existence of a universal constant $a 
> 0$ such that if $k = \frac{1}{2} n \log n + cn$ with $c > 0$, then
\num \begin{equation} \label{2.7.2}
\mbox{$\frac{1}{4}$} \sum_{[\lambda]} d_{[\lambda]}^2 \left[ \frac{1}{n} + \frac{n-1}{n} r(\lambda) \right]^{2k} \ \ \leq \ \
a^2 e^{-4c}.
\end{equation}

\noindent
This completes the proof of Theorem~\ref{2.1.3}. \vspace{.5pc}

Theorem~\ref{2.1.3} shows that $k = \frac{1}{2} n \log n + cn$ steps are sufficient for the (normalized) $\ell^2$ distance, and
hence the total variation distance, to become small.  That $k = \frac{1}{2} n \log n - cn$ steps are also \emph{necessary} is a
result of the following, which is a solution to Exercise 13 in Section D of Chapter 3 of Diaconis (1988).  The method employed in
the proof is a standard technique used in problems of this sort.

\begin{theorem} \label{2.7.3}
Let $P$ and $U$ be the probability measures on the symmetric group $S_n$ defined in (\ref{2.1.1}) and (\ref{2.1.2}),
respectively.  Let $k = \frac{1}{2} n \log n - cn$ be a nonnegative integer, with $c > 0$.  Then there exists a universal constant
$\tilde{a} > 0$ such that

\[ \mbox{$\frac{1}{2}$} \left( n! \right)^{1/2} \| P^{*k} - U \|_2 \ \ \geq \ \ \| P^{*k} - U \|_{\mbox{\rm \scriptsize TV}} \
\ \geq \ \ 1 - \tilde{a}e^{-2c}. \]

\end{theorem}

\proof{Proof} Let $\chi$ be the character of the representation $\rho_{[n-1,1]}$ of $S_n$; this representation corresponds to the
largest term in the summation from the proof of Theorem~\ref{2.1.3}.  Since any element of $S_n$ is conjugate to its inverse, it
follows from parts (b) and (c) of Lemma~\ref{2.4.1} that $\chi$ is real.

Under the uniform measure $U$, it follows from Proposition~\ref{2.4.4} that

\[ E_U( \chi ) \ = \ \frac{1}{n!} \sum_{\pi \in S_n} \chi(\pi) \ = \ \langle \chi, \chi_0 \rangle_{S_n} \ = \ 0, \]

\noindent
where $\chi_0$ is the character of the trivial representation $\rho_{[n]}$, and that

\[ \mbox{Var}_U( \chi ) \ = \ E_U( \chi^2 ) \ = \ \frac{1}{n!} \sum_{\pi \in S_n} \chi(\pi)^2 \ = \ \langle \chi, \chi
\rangle_{S_n} \ = \ 1. \]

It follows from Lemma~\ref{2.7.1} that $d_{[n-1,1]} = n-1$ and $r([n-1,1]) = \frac{n-3}{n-1}$.  Thus under the $k$-fold
convolution measure $P^{*k}$, it follows from Lemma~\ref{2.5.1} and (the calculations in the proof in Section~\ref{2.7} of)
Theorem~\ref{2.1.3} that

\[ E_{P^{*k}}( \chi ) \ = \ \sum_{\pi \in S_n} P^{*k}(\pi) \chi(\pi) \ = \ \mbox{\rm tr} \displaystyle \sum_{\pi \in S_n}
P^{*k}(\pi) \rho(\pi) \ = \ \mbox{\rm tr} \ \displaystyle \widehat{P^{*k}}(\rho) \ = \ (n-1) \left( 1 - \mbox{$\frac{2}{n}$}
\right)^k, \]

\noindent
where $\rho = \rho_{[n-1,1]}$.

In order to determine Var$_{P^{*k}}( \chi )$, we must now calculate $E_{P^{*k}}( \chi^2 )$.  It follows from Lemma~\ref{2.4.2} that
$\chi^2$ is the character of the representation $\rho_{[n-1,1]} \otimes \rho_{[n-1,1]}$.  Recall from Proposition~\ref{2.3.4} that
every representation is the direct sum of irreducible representations; in this case, it follows from the example following Lemma
2.9.16 of James and Kerber (1981) that we have the isomorphism

\[ \rho_{[n-1,1]} \otimes \rho_{[n-1,1]} \ \cong \ \rho_{[n]} \oplus \rho_{[n-1,1]} \oplus \rho_{[n-2,2]} \oplus \rho_{[n-2,1,1]}.
\]

\noindent
So it follows also from Lemma~\ref{2.4.2} that

\[ E_{P^{*k}}( \chi^2 ) \ = \ E_{P^{*k}}( \chi_{[n]} ) \ + \ E_{P^{*k}}( \chi_{[n-1,1]} ) \ + \ E_{P^{*k}}( \chi_{[n-2,2]} ) \ + \ 
E_{P^{*k}}( \chi_{[n-2,1,1]} ). \]

\noindent
As above, it follows from Lemma~\ref{2.7.1} that $d_{[n-2,2]} = \frac{1}{2}n(n-3)$ and $r([n-2,2]) = \frac{n-4}{n}$ and that
$d_{[n-2,1,1]} = \frac{1}{2}(n-1)(n-2)$ and $r([n-2,1,1]) = \frac{n-5}{n-1}$.  Thus under the $k$-fold convolution measure $P^{*k}$, 
it follows from Lemma~\ref{2.5.1} and (the calculations in the proof in Section~\ref{2.7} of) Theorem~\ref{2.1.3} that

\[ \begin{array}{rcl}

E_{P^{*k}}( \chi_{[n]} )  &  =  &  1,  \vspace{1pc} \\

E_{P^{*k}}( \chi_{[n-2,2]} )  &  =  &  \mbox{$\frac{1}{2}$} n(n-3) \left( 1 - \mbox{$\frac{2}{n}$} \right)^{2k}, \ \mathrm{and} 
\vspace{1pc} \\

E_{P^{*k}}( \chi_{[n-2,1,1]} )  &  =  &  \mbox{$\frac{1}{2}$} (n-1)(n-2) \left( 1 - \mbox{$\frac{4}{n}$} \right)^k.

\end{array} \]

\noindent
These results combine to show that

\[ \begin{array}{rcl}

\mbox{Var}_{P^{*k}}( \chi )  &  =  &  1 \ + \ (n-1) \left( 1 - \mbox{$\frac{2}{n}$} \right)^k \ + \ \mbox{$\frac{1}{2}$}
(n-1)(n-2) \left( 1 - \mbox{$\frac{4}{n}$} \right)^k  \vspace{1pc} \\

&  &  - \ \ \mbox{$\frac{1}{2}$} (n^2 - n + 2) \left( 1 - \mbox{$\frac{2}{n}$} \right)^{2k}.

\end{array} \]

Let $k = \frac{1}{2} n \log n - cn$.  By elementary calculus, $x \leq -\log(1-x) \leq \frac{x}{1-x}$ for $0 \leq x < 1$.  Thus, if
$n \geq 3$ and $c \geq 0$,

\[ \begin{array}{rcl}

\displaystyle E_{P^{*k}}( \chi )  &  =  &  \displaystyle (n-1) \left( 1 - \mbox{$\frac{2}{n}$} \right)^k \ \ \geq \ \ (n-1)
e^{-2k/(n-2)} \vspace{1pc} \\

&  =  &  \displaystyle \left( 1 - \mbox{$\frac{1}{n}$} \right) \left( \mbox{$\frac{1}{n}$} \right)^{2/(n-2)} e^{2c} e^{4c/(n-2)} \ \
\geq \ \ \mbox{$\frac{2}{27}$} e^{2c}.

\end{array} \]

\noindent
where we note that, for $n \geq 3$, $\left( 1 - \mbox{$\frac{1}{n}$} \right) \left( \frac{1}{n} \right)^{2/(n-2)}$ is increasing and
$e^{4c/(n-2)} \geq 1$.

In order to bound the variance, notice that

\[ \begin{array}{rcl}

\displaystyle \mbox{$\frac{1}{2}$} (n-1)(n-2) \left( 1 - \mbox{$\frac{4}{n}$} \right)^k  &  \leq  & \displaystyle
\mbox{$\frac{1}{2}$} (n-1)(n-2) e^{-4k/n} \ = \ \mbox{$\frac{1}{2}$} \left( 1 - \mbox{$\frac{3}{n}$} + \mbox{$\frac{2}{n^2}$}
\right) e^{4c} \ \ \ \mathrm{and} \vspace{1pc} \\

\displaystyle \mbox{$\frac{1}{2}$} (n^2 - n + 2) \left( 1 - \mbox{$\frac{2}{n}$} \right)^{2k}  &  \geq  & \displaystyle
\mbox{$\frac{1}{2}$} (n^2 - n + 2) e^{-4k/(n-2)} \vspace{1pc} \\ 

&  =  &  \displaystyle \mbox{$\frac{1}{2}$} \left( 1 - \mbox{$\frac{1}{n}$} + \mbox{$\frac{2}{n^2}$} \right) \left(
\mbox{$\frac{1}{n}$} \right)^{4/(n-2)} e^{4c} e^{8c/(n-2)}.

\end{array} \]

\noindent
Thus

\[ \mbox{$\frac{1}{2}$} (n-1)(n-2) \left( 1 - \mbox{$\frac{4}{n}$} \right)^k \ - \ \mbox{$\frac{1}{2}$} (n^2 - n + 2) \left( 1 -
\mbox{$\frac{2}{n}$} \right)^{2k} \ \leq \ 0 \]

\noindent
when

\[ \mbox{$\frac{1}{2}$} \left( 1 - \mbox{$\frac{3}{n}$} + \mbox{$\frac{2}{n^2}$} \right) e^{4c} \ - \ \mbox{$\frac{1}{2}$} \left( 1
- \mbox{$\frac{1}{n}$} + \mbox{$\frac{2}{n^2}$} \right) \left( \mbox{$\frac{1}{n}$} \right)^{4/(n-2)} e^{4c} e^{8c/(n-2)} \ \leq \
0. \]

\noindent
This restriction is equivalent to 

\[ c \ \leq \ \mbox{$\frac{1}{2}$} \log n \ - \ \left( \mbox{$\frac{n-2}{8}$} \right) \log \left( 1 - \mbox{$\frac{3}{n}$} +
\mbox{$\frac{2}{n^2}$} \right) \ + \ \left( \mbox{$\frac{n-2}{8}$} \right) \log \left( 1 - \mbox{$\frac{1}{n}$} +
\mbox{$\frac{2}{n^2}$} \right). \]

\noindent
Thus, when $c \in (0, \mbox{$\frac{1}{2}$} \log n]$,

\[ \mbox{Var}_{P^{*k}}( \chi ) \ \leq \ 1 \ + \ (n-1) \left( 1 - \mbox{$\frac{2}{n}$} \right)^k \ \leq \ 1 \ + (n-1) e^{-2k/n} \ =
\ 1 + \ \left( 1 - \mbox{$\frac{1}{n}$} \right) e^{2c} \ \leq \ 1 + \ e^{2c}. \]

Now define $A_{\alpha} := \{ \pi \in S_n : |\chi(\pi)| \leq \alpha \}$.  It follows from Chebyshev's inequality that $\displaystyle
U(A_{\alpha}) \geq 1 - \frac{1}{\alpha^2}$ and that $\displaystyle P^{*k}(A_{\alpha}) \leq \frac{1 + e^{2c}}{(\mbox{$\frac{2}{27}$}
e^{2c} - \alpha)^2}$, provided $0 \leq \alpha < \mbox{$\frac{2}{27}$} e^{2c}$.  Then

\[ \mbox{$\frac{1}{2}$} \left( n! \right)^{1/2} \| P^{*k} - U \|_2 \ \ \geq \ \ \| P^{*k} - U \|_{\mbox{\rm \scriptsize TV}} \
\ \geq \ \ 1 \ - \ \frac{1}{\alpha^2} \ - \ \frac{1 + e^{2c}}{(\mbox{$\frac{2}{27}$} e^{2c} - \alpha)^2}. \]

\noindent
Choosing $\alpha = \frac{1}{27} e^{2c}$ shows that

\[ \mbox{$\frac{1}{2}$} \left( n! \right)^{1/2} \| P^{*k} - U \|_2 \ \ \geq \ \ \| P^{*k} - U \|_{\mbox{\rm \scriptsize TV}} \
\ \geq \ \ 1 \ - \ 729 e^{-2c} \ - \ 1458 e^{-4c} \ \ \geq \ \ 1 \ - \ 2187 e^{-2c}, \]

\noindent
which completes the proof. \qed

The upper bound in Theorem~\ref{2.1.3}, taken together with the lower bound in Theorem~\ref{2.7.3}, gives an example of the 
so-called  ``cutoff phenomenon.''  The total variation distance after $k$ steps is nearly 1 until $k$ reaches about $\frac{1}{2} n
\log n$ and then drops precipitously toward 0, the dropoff occurring on the relatively small scale of $n$.  For further discussion
of the cutoff phenomenon see Diaconis (1988).

\newpage

\section{Random Walk on the Complete Monomial Groups.}  \label{3}

\subsection{Introduction.} \label{3.1}

We now extend the idea of random transpositions of $n$ cards, introduced in Section~\ref{2.1}, to a set of $n$ decks of $m$ cards each
and beyond.  Imagine $n$ decks of cards, labeled 1 through $n$, in sequential order, each with its $m$ cards in sequential order.
Independently choose two integers $p$ and $q$, uniformly from $\{1, 2, \ldots, n\}$.   

If $p \neq q$, transpose the decks in positions $p$ and $q$.  Then, independently of the choice of $p$ and $q$ and uniformly (i.e.,
with probability $\frac{1}{|G|} = \frac{1}{m!}$ each), permute the deck terminating in position $p$ by a permutation in $G = S_m$;
and independently, also uniformly, permute the deck terminating in position $q$.  This procedure is denoted by $(\vec{v};\tau)$,
where $\tau \in S_n$ is the transposition $(p\ q)$ and the only possible non-identity entries of $\vec{v} \in G^n = S_m^n$ are in
positions $p$ and $q$.  The element of $\vec{v}$ in position $p$ (resp., $q$) is $\pi \in G = S_m$ if the deck terminating in
position $p$ (resp., $q$) is permuted by $\pi \in G = S_m$.

If $p = q$ (which occurs with probability $1/n$), leave the decks in their current positions.  Then, again independently and 
uniformly, permute the deck in position $p = q$ by a permutation in $G = S_m$.  If the order of the deck is changed, this action is
denoted by $(\vec{u};e)$, where $e \in S_n$ is the identity permutation and the only non-identity entry of $\vec{u} \in G^n = S_m^n$
is in position $p = q$.  If the order of the deck is not changed, then this action is of course the identity, which is denoted by
$(\vec{e};e)$. 

If this process is repeated many times, the decks will appear to be in random order, and each deck will appear to be randomly
permuted.  

The specific example above was for motivational purposes only.  In this section we will actually examine a random walk on $G \wr
S_n$ for any group $G$, not just the symmetric group $S_m$.  In this more general setting, the example above is equivalent to
beginning with a vector $(e, \ldots, e) \in G^n$, where each $e \in G$ is the identity element.  Two elements of this vector are then
transposed as the decks of cards were above.  The transposed elements of this vector are then multiplied by elements of $G$ as the
individual decks were permuted above.

We refer to the process on $G~\wr~S_n$ described above as the \emph{independent shuffles} random walk, retaining use of the word
``shuffles'' even when $G$ is not necessarily $S_m$.  A similar process, known as the \emph{paired shuffles} random walk, will be
introduced in Section~\ref{3.7}.

In the special case of the generalized symmetric group $\mathbb{Z}_m~\wr~S_n$, imagine $n$ wheels, each of which may stop at one of
$m$ values.  The process described above randomly transposes two of the $n$ wheels and then independently spins the transposed wheels.
We thus refer to the process in this case as the \emph{independent spins} random walk.

In the special case of the hyperoctahedral group $\mathbb{Z}_2~\wr~S_n$, imagine $n$ cards, each with an orientation (up or down).
The process described above randomly transposes two of the $n$ cards and then independently flips the transposed cards.  We thus refer
to the process in this case as the \emph{independent flips} random walk.  Schoolfield (1998) used the comparison technique to analyze 
two random walks by comparing them to the independent flips random walk.

The independent shuffles random walk may be modeled formally by a probability measure $P$ on the complete monomial group $G \wr
S_n$.  Since there are ${n \choose 2}$ transpositions of $n$ elements, then there are $|G|^2 \cdot {n \choose 2} =
\frac{1}{2}|G|^2n(n-1)$ elements of the form $(\vec{v};\tau)$.  There are also $n(|G|-1)$ elements of the form $(\vec{u};e)$ with
$\vec{u} \neq \vec{e}$.  We may thus define the following probability measure on the set of all elements of $G~\wr~S_n$:
\num \begin{equation} \label{3.1.1}
\begin{array}{rcll}

P(\vec{e};e) & = & \displaystyle \frac{1}{|G|n}, & \vspace{.5pc} \\

P(\vec{u};e) & = & \displaystyle \frac{1}{|G|n^2} & \mbox{where $\vec{u} \neq \vec{e} \in G^n$}, \vspace{.5pc} \\

P(\vec{v};\tau) & = & \displaystyle \frac{2}{|G|^2 n^2} & \mbox{where $\vec{v} \in G^n$}, \vspace{.5pc} \\

P(\vec{x};\pi) & = & 0 & \mbox{otherwise},

\end{array}
\end{equation}

\noindent
where there is only one non-identity entry of $\vec{u} \in G^n$, and where if $\tau \in S_n$ is the transposition $(p\ q)$ then the
only possible non-identity entries of $\vec{v} \in G^n$ are in positions $p$ and $q$.  In the special case of the hyperoctahedral
group $\mathbb{Z}_2~\wr~S_n$, we refer to the elements $(\vec{u};e)$ as \emph{signed identities} and the elements $(\vec{v};\tau)$ as
\emph{signed transpositions}.

Since there are $|G|^n \cdot n!$ elements in $G~\wr~S_n$, the uniform probability measure is given by
\num \begin{equation} \label{3.1.2}
U(\vec{x};\pi) = \frac{1}{|G|^n \cdot n!} \ \ \ \mbox{for $(\vec{x};\pi) \in G~\wr~S_n$}.
\end{equation}

The following result establishes an upper bound on both the total variation distance and the $\ell^2$ distance between $P^{*k}$ and
$U$.  We establish an analogous result for the paired shuffles random walk as Theorem~\ref{3.7.3}.

\begin{theorem} \label{3.1.3}
Let $P$ and $U$ be the probability measures on the complete monomial group $G~\wr~S_n$ defined in (\ref{3.1.1}) and
(\ref{3.1.2}), respectively.  Let $k = \frac{1}{2} n \log n + \frac{1}{4} n \log(|G|-1) + cn$.  Then there exists a universal
constant $b>0$ such that

\[ \| P^{*k} - U \|_{\mbox{\rm \scriptsize TV}} \ \ \leq \ \ \mbox{$\frac{1}{2}$} \left(|G|^n n!\right)^{1/2} \| P^{*k} - U
\|_2 \ \ \leq \ \ be^{-2c} \ \ \ \mathrm{for\ all\ } \mbox{$c > 0$}. \]

\end{theorem}

In the following sections we establish the results necessary to prove this theorem and an analogous theorem for the paired
shuffles random walk.  In Section~\ref{3.2} we study the basic properties of the complete monomial groups.  In Section~\ref{3.3}
we show that the probability measure defined above is constant on conjugacy classes.  In Section~\ref{3.4} we identify all of the
irreducible representations of the complete monomial groups and in Section~\ref{3.5} we calculate the characters of these
irreducible representations.  In Section~\ref{3.6} we calculate the Fourier transform of the probability measure defined in
(\ref{3.1.1}) and prove Theorem~\ref{3.1.3} along with a matching $\ell^2$ lower bound.  In Section~\ref{3.7} we perform a similar
analysis of the paired shuffles random walk; however, the results are quite different in the case that $G$ is nonabelian.

\subsection{The Complete Monomial Groups.} \label{3.2}

The \emph{complete monomial group} $G~\wr~S_n$ is the wreath product of the group $G$ with the symmetric group $S_n$.  Special cases
include the \emph{hyperoctahedral group} $\mathbb{Z}_2~\wr~S_n$, the \emph{generalized symmetric group} $\mathbb{Z}_m~\wr~S_n$, and
$S_m~\wr~S_n$.  The elements of $G~\wr~S_n$ may be represented as $(\vec{x};\pi) \in\ G^n~\times~S_n$.  It then follows that the order
of the complete monomial group is $|G~\wr~S_n| = |G|^n \cdot n!$.   

Each element $(\vec{x};\pi) \in G~\wr~S_n$ acts on a vector $\vec{w} \in G^n$ by first permuting its elements according to the
permutation $\pi$ and then left-multiplying the elements of the permuted vector by the elements of $\vec{x}$, entry by entry; i.e.,
for any $(\vec{x};\pi) \in G~\wr~S_n$ and any $\vec{w} \in G^n$,

\[ (\vec{x};\pi)(\vec{w}) = \left( x_1 \cdot w_{\pi^{-1}(1)}, \ldots, x_n \cdot w_{\pi^{-1}(n)} \right) \in G^n. \]

\noindent
In keeping with the decks of cards analogy, each element $(\vec{x};\pi) \in\ S_m~\wr~S_n$ first permutes the $n$ decks of cards
according to $\pi$ and then permutes each deck of cards according to the entry of $\vec{x}$ at its new index.   

Considering the actions, described above, of the elements of the complete monomial group on vectors in $G^n$, it follows that the
product of two elements  $(\vec{x};\pi), (\vec{y};\sigma) \in G~\wr~S_n$ is given by

\[ (\vec{y};\sigma) \cdot (\vec{x};\pi) = \left( y_1 \cdot x_{\sigma^{-1}(1)}, \ldots, y_n \cdot x_{\sigma^{-1}(n)}; \sigma \pi
\right). \]

\noindent
Thus the identity element is $\left( e, \ldots, e; e \right)$ where $(e, \ldots, e) \in G^n$ and $e$ is the identity permutation
in $S_n$.  It also follows that the inverse of any element $(\vec{x};\pi) \in G~\wr~S_n$ is 

\[ (\vec{x};\pi)^{-1} = \left( x_{\pi(1)}^{-1}, \ldots, x_{\pi(n)}^{-1}; \pi^{-1} \right). \]

There are two very important subgroups of $G~\wr~S_n$ which should be considered.  The subgroup consisting of elements of the form
$\{ (\vec{x};e) : \vec{x} \in G^n \}$ is isomorphic to $G^n$ and is a normal subgroup of $G~\wr~S_n$. The subgroup consisting of
elements of the form $\{ (e, \ldots, e; \pi) : \pi \in S_n \}$ is isomorphic to $S_n$.  Notice that any element of $G~\wr~S_n$ can
be written uniquely as the product of a single element from each of these two subgroups.   

\subsection{Conjugacy Classes.} \label{3.3}

Recall that any permutation $\pi \in S_n$ can be written as the product of at most $n$ disjoint cyclic factors.  (For any index
not moved by $\pi$, there is considered to be a cyclic factor of length one.)  For any $\pi \in S_n$ there is an associated
$n$-dimensional vector, called the \emph{cycle type} of $\pi$, which lists the number of disjoint cyclic factors of $\pi$ of each
possible length. 

For any $(\vec{x};\pi) \in G~\wr~S_n$, the permutation $\pi$ can be written as the product of disjoint cyclic factors, as described
above.  Suppose that $\pi$ is the product of $j$ disjoint cyclic factors of lengths $\{k_1, \ldots, k_j\}$, respectively, with $k_1 + 
\cdots + k_j = n$.  Then $\pi$ can be written as 

\[ \pi = (i_1^{(1)}\ i_2^{(1)}\ \cdots\ i_{k_1}^{(1)})\ (i_1^{(2)}\ i_2^{(2)}\ \cdots\ i_{k_2}^{(2)})\ \cdots\ (i_1^{(j)}\ 
i_2^{(j)}\ \cdots\ i_{k_j}^{(j)}) \]

\noindent
where each $i_m^{(\ell)} \in \{1, 2, \ldots, n\}$.  For each disjoint cyclic factor of $\pi$, we may define

\[ g_{\ell}(\vec{x};\pi) = x_{i_1^{(\ell)}} \cdot x_{i_2^{(\ell)}} \cdot \, \cdots \, \cdot x_{i_{k_{\ell}}^{(\ell)}}, \] 

\noindent
which is called the $\ell$th \emph{cycle product} of $(\vec{x};\pi)$.  Notice that $g_{\ell}(\vec{x};\pi) \in G$ for all $1 \leq
\ell \leq j$.   

Suppose that $G$ has $s$ conjugacy classes $C_1, C_2, \ldots, C_s$, with $C_1$ being the conjugacy class of the identity element $e
\in G$.  By calculating $g_{\ell}(\vec{x};\pi)$ for all $1 \leq \ell \leq j$, we may construct an $s~\times~n$ \emph{type matrix}
$a(\vec{x};\pi)$ for $(\vec{x};\pi)$ in the following manner.  Let the $k$th entry of the $j$th row be the number of cyclic factors
of length $k$ contained in $\pi$ for which $g_{\ell}(\vec{x};\pi) \in C_j$.  According to Theorem 4.2.8 of James and Kerber (1981),
two elements of $G~\wr~S_n$ are conjugate if and only if their type matrices are identical.  Notice that the vector of column sums
of $a(\vec{x};\pi)$ is the cycle type vector of $\pi$.   

We are primarily interested in the type matrices of the elements in the support of the probability measure defined in (\ref{3.1.1}).
In the identity element $(\vec{e};e) \in\ G~\wr~S_n$, the identity permutation $e \in S_n$ is the product of $n$ disjoint cyclic
factors each of length one, and each of these cyclic factors corresponds to an identity element in $\vec{e} = (e, \ldots, e) \in
G^n$.  Thus the type matrix of the identity element is the $s~\times~n$ matrix with $a_{11}(\vec{e};e) = n$ and the remaining
entries all zeros.   

In the special case of the hyperoctahedral group $\mathbb{Z}_2~\wr~S_n$, we have $s = 2$.  Thus the type matrix of the identity
element is 

\onespace

\[ a(\vec{0};e) = \left( \begin{array}{c c c c c} n & 0 & 0 & \cdots & 0 \\ 0 & 0 & 0 & \cdots & 0 \end{array} \right). \]

\twospace

Notice that in each $(\vec{u};e) \in G~\wr~S_n$ with $\vec{u} \neq \vec{e}$, the identity permutation $e \in S_n$ is again the
product of $n$ disjoint cyclic factors each of length one.  However, $\vec{u} \in G^n$ has a single non-identity entry, which
corresponds to one of the cyclic factors of the identity permutation.  Thus, there are $s-1$ different type matrices for these
elements, each with $a_{11}(\vec{u};e) = n-1$, with $a_{k1}(\vec{u};e) = 1$ if the non-identity entry of $\vec{u}$ is in $C_k$ and
the remaining entries all zeros.  Hence, the elements of the form $(\vec{u};e) \in G~\wr~S_n$ split into $s-1$ conjugacy classes.

In the special case of the hyperoctahedral group $\mathbb{Z}_2~\wr~S_n$, the type matrix of each of the signed identities, which
together form a single conjugacy class, is

\onespace

\[ a(\vec{u};e) = \left( \begin{array}{c c c c c} n-1 & 0 & 0 & \cdots & 0 \\ 1 & 0 & 0 & \cdots & 0 \end{array} \right). \]

\twospace

In each $(\vec{v};\tau) \in G~\wr~S_n$, the transposition $\tau \in S_n$ is the product of $n-2$ cyclic factors of length
one and one cyclic factor of length two.  Also, each of the $n-2$ cyclic factors of length one corresponds to an identity $e \in G$ 
in $\vec{v} \in G^n$, while the only possible non-identity entries in $\vec{v} \in G^n$ correspond to the one cyclic factor of
length two.  Thus $g_{\ell}(\vec{v};\tau) = 0$ for each of the $n-2$ cyclic factors of length one, but $g_{\ell}(\vec{v};\tau)$ can
be an arbitrary element of $G$ for the one cyclic factor of length two.  Thus, there are $s$ different type matrices for these
elements, each with $a_{11}(\vec{v};\tau) = n-2$, with $a_{k2}(\vec{v};\tau) = 1$ if $g_{\ell}(\vec{v};\tau) \in C_k$ for the one
cyclic factor of length two, and with the remaining entries all zeros.  Hence, the elements of the form $(\vec{v};\tau) \in G~\wr~S_n$
split into $s$ conjugacy classes.

In the special case of the hyperoctahedral group $\mathbb{Z}_2~\wr~S_n$, the calculation of the type matrix splits the signed
transpositions into two sets.  Let us refer to a signed transposition which flips neither or both of the cards as an \emph{even}
transposition and designate it as $(\vec{v};\tau^+)$.  We will refer to one which flips either, but not both, of the cards as an
\emph{odd} transposition and designate it as $(\vec{v};\tau^-)$.  Thus the type matrices of the even transpositions and the odd
transpositions are, respectively,

\onespace

\[ a(\vec{v};\tau^+) = \left( \begin{array}{c c c c c} n-2 & 1 & 0 & \cdots & 0 \\ 0 & 0 & 0 & \cdots & 0 \end{array} \right),
\ \ \ \ \ 
   a(\vec{v};\tau^-) = \left( \begin{array}{c c c c c} n-2 & 0 & 0 & \cdots & 0 \\ 0 & 1 & 0 & \cdots & 0 \end{array} \right).
\]

\twospace

We have now identified the type matrices (and, hence, the conjugacy classes) of all the elements in the support of the probability
measure defined in (\ref{3.1.1}).  Furthermore, we have established the following.

\begin{lemma} \label{3.3.1}
The probability measure defined in (\ref{3.1.1}) is constant on conjugacy classes.   
\end{lemma}

Knowing that there are the same number of conjugacy classes as there are type matrices allows us to calculate the number of
conjugacy classes.  According to Lemma 4.2.9 of James and Kerber (1981), the number of conjugacy classes of $G~\wr~S_n$ is equal
to

\[ \sum_{(n_1, n_2, \ldots, n_s)} p(n_1) \cdot p(n_2) \cdots p(n_s) \]

\noindent
where $p(k)$ is the number of partitions of $k$ (with $p(0) := 1$) and the sum is taken over all $s$-dimensional vectors such that
$n_1 + n_2 + \cdots + n_s = n$ and $n_j \geq 0$ for all $1 \leq j \leq s$.  

We may also calculate the order of each of the conjugacy classes.  According to Lemma 4.2.10 of James and Kerber (1981), the number
of elements of a particular type $a(\vec{x};\pi)$ in $G~\wr~S_n$, and hence the number of elements of a particular conjugacy class,
is equal to

\[ \frac{|G|^n \cdot n!}{\prod_{i,j} \left[ (j|G|/|C_i|)^{a_{ij}} \cdot a_{ij}! \right]} \]

\noindent
where $a_{ij}$ is the $(i,j)$th element of $a(\vec{x};\pi)$ and $|C_i|$ is the order of the $i$th conjugacy class of $G$.  By
applying this result to the elements of $G~\wr~S_n$ whose type matrices were determined above, we have the following otherwise
obvious corollary.

\begin{corollary} \label{3.3.2}
There are $n|C_k|$ elements of the form $(\vec{u};e)$, for $2 \leq k \leq s$, and there are $\frac{1}{2}n(n-1) \cdot |G| \cdot
|C_k|$ elements of the form $(\vec{v};\tau)$, for $1 \leq k \leq s$.   
\end{corollary}

\subsection{Irreducible Representations.} \label{3.4}

We now construct a collection of irreducible representations of the complete monomial groups from irreducible representations
of certain of their subgroups.  We will later see that \emph{every} irreducible representation of $G~\wr~S_n$ can be constructed in
such a manner.  The method used is from Section 4.3 of James and Kerber (1981).  Other methods described in Section 8.2 of Serre
(1977) and in Chapter V of Simon (1996) could be used in the special case when $G$ is abelian.   

Let $G^*$ be the subgroup of $G~\wr~S_n$ consisting of all elements of the form $\{(\vec{x};e) : \vec{x} \in G^n \}$ and let $H^*$
be the subgroup of $G~\wr~S_n$ consisting of all elements of the form $\{(e, \ldots, e;\pi) : \pi \in S_n \}$.  Recall from
Section~\ref{3.2} that $G^*$ is a normal subgroup of $G~\wr~S_n$ and that $G^* \cong G^n$ and $H^* \cong S_n$.  Recall, furthermore,
that each element of $G~\wr~S_n$ can be written uniquely as a product $gh$ with $g \in G^*$ and $h \in H^*$.   

Suppose that $G$ has $s$ conjugacy classes.  Then it follows from Proposition \ref{2.3.1} that there are $s$ irreducible
representations of $G$, namely, $\rho_1, \rho_2, \ldots, \rho_s$, where we choose our labelling so that $\rho_1$ is the trivial
representation.  So it follows from Proposition \ref{2.3.5} that any irreducible representation of $G^n$ is isomorphic to an $n$-fold
tensor product 

\[ \rho^{(1)} \otimes \rho^{(2)} \otimes \cdots \otimes \rho^{(n)} \]

\noindent
where each $\rho^{(i)} \in \{\rho_1, \rho_2, \ldots, \rho_s \}$.  Notice that any irreducible representation of $G^n$ can be 
represented by a vector $\vec{x} \in\ \{1, \ldots, s\}^n$, where the entries of $\vec{x}$ correspond to the indices of the factors
in the tensor product above.   

Suppose that the tensor product forming a particular irreducible representation of $G^n$ contains $n_j$ factors of the kind $\rho_j$,
for all $1 \leq j \leq s$.  The vector $(n) = (n_1, n_2, \ldots, n_s)$ is called the $type$ of this representation.  Notice that
for each type $(n) = (n_1, n_2, \ldots, n_s)$, there is a unique representative whose first $n_1$ factors in the tensor product
are $\rho_1$, whose next $n_2$ factors are $\rho_2$, etc.  We will refer to this canonical representative as $\rho_{(n)}$.   

Since $G^n \cong G^* = \{(\vec{x};e) : \vec{x} \in G^n \}$, any representation $\rho$ of $G^n$ extends easily to a representation
of $G^*$ by setting $\rho(\vec{x};e) \equiv \rho(\vec{x})$.  Thus any representation of $G^n$ is also a representation of $G^*$.
Therefore, the collection $\left\{ \rho_{(n)} \right\}$ is a complete system of irreducible representations of $G^*$ which are of
different types.   

We now turn our attention to $H^* = \{(e, \ldots, e;\pi) : \pi \in S_n \} \cong S_n$.  For each type $(n) = (n_1, n_2,
\ldots, n_s)$, define $S_{(n)}$ to be the subgroup of $S_n$ which permutes the first $n_1$ indices among themselves, the next $n_2$
indices among themselves, etc.; but does not commingle these $s$ sets of indices.  Thus $S_{(n)} \cong S_{n_1} \times S_{n_2} \times
\cdots \times S_{n_s}$, where we define $S_{n_j} := S_1$ if $n_j = 0$ for any $1 \leq j \leq s$. 

Recall that there is a one-to-one correspondence between irreducible representations of $S_n$ and partitions $[\lambda]$ of $n$, 
where $[\lambda] = [\lambda_1, \lambda_2, \ldots, \lambda_k]$ with $\lambda_1 \geq \lambda_2 \geq \cdots \geq \lambda_k > 0$ and
$\lambda_1 + \lambda_2 + \cdots + \lambda_k = n$.  Thus any irreducible representation of $S_n$ may be denoted as $\rho_{[\lambda]}$,
where $[\lambda]$ is the corresponding partition of $n$.

It follows from Proposition \ref{2.3.5} that any irreducible representation of $S_{n_1} \times S_{n_2} \times \cdots \times S_{n_s}$ is
isomorphic to the $s$-fold tensor product

\[ \rho_{[\lambda_1]} \otimes \rho_{[\lambda_2]} \otimes \cdots \otimes \rho_{[\lambda_s]} \]

\noindent
where $(\lambda) = ([\lambda_1], [\lambda_2], \ldots, [\lambda_s])$ are partitions of $\{ n_1, n_2, \ldots, n_s \}$ and
$\{ \rho_{[\lambda_1]}, \rho_{[\lambda_2]}, \ldots,  \rho_{[\lambda_s]} \}$ are irreducible representations of $\{S_{n_1},
S_{n_2}, \ldots, S_{n_s} \}$, respectively.  We will denote this tensor product as $\rho_{(\lambda)}$.  Thus, there is a
one-to-one correspondence between irreducible representations of $S_{n_1} \times S_{n_2} \times \cdots \times S_{n_s}$ and
ordered $s$-tuples of partitions of $(n) = (n_1, n_2, \ldots, n_s)$.  Since $S_{n_1} \times S_{n_2} \times \cdots \times S_{n_s}
\cong S_{(n)}$, then any representation of $S_{n_1} \times S_{n_2} \times \cdots \times S_{n_s}$ is also a representation of
$S_{(n)}$.   

Let us define $S_{(n)}^* = \{(e, \ldots, e; \pi) : \pi \in S_{(n)} \}$.  Since $S_{(n)}$ is a subgroup of $S_n$, it follows that
$S_{(n)}^*$ is a subgroup of $G^n~\wr~S_n$.  Since $S_{(n)}^* \cong S_{(n)}$, any representation $\rho$ of $S_{(n)}$ extends easily
to a representation of $S_{(n)}^*$ by setting $\rho(e, \ldots, e; \pi) \equiv \rho(\pi)$.  Thus any representation of $S_{(n)}$ is
also a representation of $S_{(n)}^*$.  Therefore, for fixed type $(n) = (n_1, n_2, \ldots, n_s)$, the collection $\left\{
\rho_{(\lambda)} \right\}$ is a complete system of irreducible representations of $S_{(n)}^*$, when $(\lambda)$ ranges over all
ordered $s$-tuples of partitions of $(n_1, n_2, \ldots, n_s)$, respectively.   

Now consider the wreath product $G~\wr~S_{(n)}$, which is a subgroup of $G~\wr~S_n$.  The representations $\rho_{(n)} : G^*
\longrightarrow \mathrm{\textbf{GL}}(V_1)$ and $\rho_{(\lambda)} : S_{(n)}^* \longrightarrow \mathrm{\textbf{GL}}(V_2)$, where
$(\lambda) = ([\lambda_1], [\lambda_2], \ldots, [\lambda_s])$ are partitions of $(n) = (n_1, n_2, \ldots, n_s)$, respectively, may
be combined to form representations $\rho_{(n)} \otimes \rho_{(\lambda)}$ of $G~\wr~S_{(n)}$ by defining (with an abuse of the
notation ``$\otimes$,'' since this is not quite the same notion of tensor product as in Section \ref{2.3})

\[ \left( \rho_{(n)} \otimes \rho_{(\lambda)} \right) (\vec{x};\pi) (v_1 \otimes v_2) \ := \ \rho_{(n)} (\vec{x};e) (v_1) \otimes
\rho_{(\lambda)} (\vec{0};\pi) (v_2) \]

\noindent
for all $(\vec{x};\pi) \in G~\wr~S_{(n)}$ and all $v_1 \in V_1$ and $v_2 \in V_2$.  These representations $\left\{ \rho_{(n)}
\otimes \rho_{(\lambda)} \right\}$ of $G~\wr~S_{(n)}$ may now be used to induce representations of $G~\wr~S_n$.  More importantly,
we have the following, which is Theorem 4.4.3 of James and Kerber (1981).

\begin{lemma} \label{3.4.1}
The collection of representations of $G~\wr~S_n$ induced by the collection $\left\{ \rho_{(n)} \otimes \rho_{(\lambda)} \right\}$ of 
representations of $G~\wr~S_{(n)}$ is a complete collection of pairwise inequivalent and irreducible representations of $G~\wr~S_n$
if $(n) = (n_1, n_2, \ldots, n_s)$ ranges over all different types and, for fixed type $(n)$, $(\lambda) = ([\lambda_1],
[\lambda_2], \ldots, [\lambda_s])$ ranges over all ordered $s$-tuples of partitions of $(n_1, n_2, \ldots, n_s)$, respectively.  
\end{lemma}

Since the number of irreducible representations equals the number of conjugacy classes (Proposition \ref{2.3.1}), it follows from the
results in Section~\ref{3.3} that the number of irreducible representations of $G~\wr~S_n$ is equal to

\[ \sum_{(n_1, n_2, \ldots, n_s)} p(n_1) \cdot p(n_2) \cdots p(n_s), \]

\noindent
where $p(k)$ is the number of partitions of $k$ (with $p(0) := 1$) and the sum is taken over all $s$-dimensional vectors such that
$n_1 + n_2 + \cdots + n_s = n$ and $n_j \geq 0$ for all $1 \leq j \leq s$.  This is Corollary 4.4.4 of James and Kerber (1981) and 
is consistent with the results found above.   

\subsection{Irreducible Characters.} \label{3.5}

For the elements in the support of the probability measure defined in (\ref{3.1.1}), we now determine the characters of the
irreducible representations of the complete monomial group $G~\wr~S_n$ induced by the irreducible representations of $G~\wr~S_{(n)}$
found in Lemma~\ref{3.4.1}.  We do this with the aid of Lemma \ref{2.4.3}.  The Frobenius character formula, which is Theorem V.4.1 of
Simon (1996) or Theorem 12 of Serre (1977), may also be used.   

For notational purposes, let $C_k^{(\vec{u};e)}$ (with $2 \leq k \leq s$) be the conjugacy classes of $(\vec{u};e)$ in $G~\wr~S_n$,
with $k$ chosen so that the single non-identity entry of $\vec{u}$ is in conjugacy class $C_k$ of $G$, and let
$C_k^{(\vec{v};\tau)}$ (with $1 \leq k \leq s$) be the conjugacy classes of $(\vec{v};\tau)$ in $G~\wr~S_n$, with $k$ chosen so that
the product (calculated as in the determination of the type matrices in Section~\ref{3.3}) of the two possible non-identity entries
of $\vec{v}$ is in conjugacy class $C_k$ of $G$.  

\begin{lemma} \label{3.5.1}
For the elements in the support of the probability measure defined in (\ref{3.1.1}), the character of the irreducible
representation $\rho$ of the complete monomial group $G~\wr~S_n$ induced by the irreducible representation $\rho_{(n)} \otimes
\rho_{(\lambda)}$ of $G~\wr~S_{(n)}$ is given by

\[ \begin{array}{rcl}

\displaystyle \chi_{\rho}(\vec{e};e) & = & \displaystyle {n \choose n_1, \ldots, n_s} d_{\rho_1}^{n_1} \cdots d_{\rho_s}^{n_s} \cdot 
d_{[\lambda_1]} \cdots d_{[\lambda_s]} \ = \ d_{\rho}, \vspace{1pc} \\  

\displaystyle \chi_{\rho}(\vec{u};e) & = & \displaystyle d_{\rho} \sum_{j=1}^s \left( \frac{n_j}{n} \right) \left[ 
\frac{\chi_{\rho_j}(g_k)}{d_{\rho_j}} \right] \ \ \ \ \ \mathrm{for\ } \mbox{$(\vec{u};e) \in C_k^{(\vec{u};e)}$}, \vspace{1pc} \\  

\displaystyle \chi_{\rho}(\vec{v};\tau) & = & \displaystyle d_{\rho} \sum_{j=1}^s \left[ \frac{n_j(n_j-1)}{n(n-1)} \right] \cdot
\left[ \frac{\chi_{\rho_j}(g_k)}{d_{\rho_j}} \right] \cdot r(\lambda_j) \ \ \ \ \ \mathrm{for\ } \mbox{$(\vec{v};\tau) 
\in C_k^{(\vec{v};\tau)}$},

\end{array} \]

\noindent
where, for $1 \leq j \leq s$, $\chi_{\rho_j}$ is the character of the irreducible representation $\rho_j$ of $G$ and 
$\chi_{[\lambda_j]}$ is the character of the irreducible representation $\rho_{[\lambda_j]}$ of $S_{n_j}$, where $r(\lambda_j)
:= \chi_{[\lambda_j]}(\tau)/d_{[\lambda_j]}$ with transposition $\tau \in S_{n_j}$, and where $g_k$ is any element of the conjugacy
class $C_k$ of~$G$.

\end{lemma}

\proof{Proof} Let $\widetilde{G} := G~\wr~S_n$ and $\widetilde{H} := G~\wr~S_{(n)}$; we apply Lemma \ref{2.4.3} to these groups.
Notice that $|\widetilde{G}| = |G~\wr~S_n| = |G|^n \cdot n!$ and that $|\widetilde{H}| = |G~\wr~S_{(n)}| = \left|G~\wr~\left(S_{n_1}
\times \cdots S_{n_s} \right) \right| = |G|^n \cdot n_1! \cdots n_s!$.  Thus

\[ \displaystyle \frac{|\widetilde{G}|}{|\widetilde{H}|} = \frac{|G|^n \cdot n!}{|G|^n \cdot n_1! \cdots n_s!} = {n \choose n_1,
\ldots, n_s}. \]

Recall, in the representation $\rho_{(n)} \otimes \rho_{(\lambda)}$ of $\widetilde{H} = G~\wr~S_{(n)}$, that $\rho_{(n)}$ is
the $n$-fold tensor product

\[ \rho^{(1)} \otimes \rho^{(2)} \otimes \cdots \otimes \rho^{(n)}, \]

\noindent
where each $\rho^{(i)} \in \{\rho_1, \rho_2, \ldots, \rho_s \}$ is an irreducible representation of $G$ and the first $n_1$ factors
are $\rho_1$, the next $n_2$ factors are $\rho_2$, etc.  Also recall that $\rho_{(\lambda)}$ is

\[ \rho_{[\lambda_1]} \otimes \rho_{[\lambda_2]} \otimes \cdots \otimes \rho_{[\lambda_s]}, \] 

\noindent
where $\rho_{[\lambda_1]}, \rho_{[\lambda_2]}, \ldots, \rho_{[\lambda_s]}$ are irreducible representations of $S_{n_1}, S_{n_2},
\ldots, S_{n_s}$, respectively.  It then follows from Lemma \ref{2.4.2} that the character of the representation $\rho_{(n)}
\otimes\rho_{(\lambda)}$ is given by

\[ \chi_{\widetilde{H}} = \chi^{(1)} \cdots \chi^{(n)} \cdot \chi_{[\lambda_1]} \cdots \chi_{[\lambda_s]}. \]

We begin with the identity element $(\vec{e};e) \in \widetilde{G} = G~\wr~S_n$.  Since the identity element forms a
singleton conjugacy class, we have $\ell = 1$.  Furthermore, since the only element of $\widetilde{H}$ to which it is conjugate in
$\widetilde{G}$ is the identity element $(\vec{e};e, \ldots, e) \in \widetilde{H} = G~\wr~S_{(n)}$, we have $t = 1$ and $k_1 = 1$.   

Recall that the character at the identity of any representation is the dimension of the representation.  Thus

\[ \chi^{(1)}(e) \cdots \chi^{(n)}(e) = d_{\rho_1}^{n_1} \cdots d_{\rho_s}^{n_s}. \]

\noindent
So the character of the identity element $(\vec{e};e, \ldots, e) \in \widetilde{H} = G~\wr~S_{(n)}$ is

\[ d_{\rho_1}^{n_1} \cdots d_{\rho_s}^{n_s} \cdot \chi_{[\lambda_1]}(e) \cdots \chi_{[\lambda_s]}(e). \]

Therefore, it follows from Lemma \ref{2.4.3} that the character at the identity element in the induced irreducible representation 
$\rho$ of $G~\wr~S_n$ is

\[ \chi_{\rho}(\vec{e};e) = {n \choose n_1, \ldots, n_s} d_{\rho_1}^{n_1} \cdots d_{\rho_s}^{n_s} \cdot \chi_{[\lambda_1]}(e) \cdots
\chi_{[\lambda_s]}(e), \] 

\noindent
and, furthermore, that the dimension of the induced irreducible representation $\rho$ of $G~\wr~S_n$ is

\[ d_{\rho} = {n \choose n_1, \ldots, n_s} d_{\rho_1}^{n_1} \cdots d_{\rho_s}^{n_s} \cdot d_{[\lambda_1]} \cdots d_{[\lambda_s]}. \] 

We now consider the elements $(\vec{u};e) \in \widetilde{G} = G~\wr~S_n$ with $\vec{u} \neq \vec{e}$, which comprise $s-1$
conjugacy classes: $C_2^{(\vec{u};e)}, C_3^{(\vec{u};e)}, \ldots, C_s^{(\vec{u};e)}$, where the indices are chosen so that, for
$C_k^{(\vec{u};e)}$, the single non-identity entry of $\vec{u}$ is in conjugacy class $C_k$ of $G$.  For a particular conjugacy
class $C_k^{(\vec{u};e)}$, since there are $n|C_k|$ elements in $\widetilde{G}$, we have $\ell = n|C_k|$.  Each of the $s-1$
conjugacy classes splits into $s$ classes in $\widetilde{H}$:  in the first class, the only non-identity entry of $\vec{u}$ is one
of its first $n_1$ entries; in the second class, the only non-identity entry of $\vec{u}$ is one of its next $n_2$ entries; etc.
Thus $t = s$ with $k_j = n_j|C_k|$ for all $1 \leq j \leq s$.   

For a particular conjugacy class $C_k$ of $G$, let $g_k$ be a representative element.  Then the characters of the $s$ classes of
elements $(\vec{u};e, \ldots, e) \in \widetilde{H} = G~\wr~S_{(n)}$ are

\[ d_{\rho_1}^{n_1} \cdots d_{\rho_s}^{n_s} \cdot d_{[\lambda_1]} \cdots d_{[\lambda_s]} \cdot 
\frac{\chi_{\rho_j}(g_k)}{d_{\rho_j}}, \]

\noindent
for $1 \leq j \leq s$.   

Therefore, it follows from Lemma \ref{2.4.3} that the character at $(\vec{u};e)$ of the induced irreducible representation $\rho$ of $G
\wr S_n$ is

\[ \chi_{\rho}(\vec{u},e) = {n \choose n_1, \ldots, n_s} \sum_{j=1}^s \frac{n_j}{n}  d_{\rho_1}^{n_1} \cdots
d_{\rho_s}^{n_s} \cdot d_{[\lambda_1]} \cdots d_{[\lambda_s]} \cdot \frac{\chi_{\rho_j}(g_k)}{d_{\rho_j}}, \]

\noindent
for $2 \leq k \leq s$, from which the desired result follows.   

Finally, we consider the elements $(\vec{v};\tau) \in \widetilde{G} = G~\wr~S_n$, which comprise $s$ conjugacy classes:
$C_1^{(\vec{v};\tau)}, C_2^{(\vec{v};\tau)}, \ldots, C_s^{(\vec{v};\tau)}$, where the indices are chosen so that, for
$C_k^{(\vec{v};\tau)}$, the product (calculated as in the determination of the type matrices in Section~\ref{3.3}) of the two
possible non-identity entries of $\vec{v}$ is in conjugacy class $C_k$ of $G$.  Since there are $\frac{1}{2}n(n-1) \cdot |G| \cdot
|C_k|$ elements in the conjugacy class $C_k^{(\vec{v};\tau)}$ of $\widetilde{G}$, we have $\ell = \frac{1}{2}n(n-1) \cdot |G| \cdot
|C_k|$.   

Each of the $s$ conjugacy classes splits into $s$ classes in $\widetilde{H}$:  in the first class, $\tau$ transposes two of the
first $n_1$ elements leaving other elements fixed; in the second class, $\tau$ transposes two of the next $n_2$ elements leaving the
other elements fixed; etc.  Also, in the first class, the only non-identity entries of $\vec{v}$ are in its first $n_1$ entries; in
the second class, the only non-identity entries of $\vec{v}$ are in its next $n_2$ entries; etc.  Thus $t = s$ with $k_j =
\frac{1}{2}n_j(n_j-1) \cdot |G| \cdot |C_k|$ for all $1 \leq j \leq s$.   

For a particular conjugacy class $C_k$ of $G$, let $g_k$ be a representative element.  Then the characters of the $s$ classes of
elements $(\vec{v};\tau) \in \widetilde{H} = G~\wr~S_{(n)}$ are

\[ d_{\rho_1}^{n_1} \cdots d_{\rho_s}^{n_s} \cdot d_{[\lambda_1]} \cdots d_{[\lambda_s]} \cdot \frac{\chi_{\rho_j}(g_k)}{d_{\rho_j}}
\cdot \frac{\chi_{[\lambda_j]}(\tau)}{d_{[\lambda_j]}}, \]

\noindent
for $1 \leq j \leq s$.   

Therefore, it follows from Lemma \ref{2.4.3} that the character at $(\vec{v};\tau)$ of the induced irreducible representation $\rho$
of $G~\wr~S_n$ is

\[ \chi_{\rho}(\vec{v},\tau) = {n \choose n_1, \ldots, n_s} \sum_{j=1}^s \frac{n_j(n_j-1)}{n(n-1)} d_{\rho_1}^{n_1} \cdots
d_{\rho_s}^{n_s} \cdot d_{[\lambda_1]} \cdots d_{[\lambda_s]} \cdot \frac{\chi_{\rho_j}(g_k)}{d_{\rho_j}} \cdot
\frac{\chi_{[\lambda_j]}(\tau)}{d_{[\lambda_j]}}, \]

\noindent
for $1 \leq k \leq s$, from which the desired result follows.  $\qed$   

\subsection{Analysis of the Independent Shuffles Random Walk.} \label{3.6}

In order to continue our analysis of the independent shuffles random walk introduced in Section~\ref{3.1}, we must now calculate
the Fourier transform of $P$ at each irreducible representation of the complete monomial group $G~\wr~S_n$.

\begin{lemma} \label{3.6.1}
Let $P$ be the probability measure on $G~\wr~S_n$ defined in (\ref{3.1.1}).  For the irreducible representation $\rho$ of
$G~\wr~S_n$ induced by the representation $\rho_{(n)} \otimes \rho_{(\lambda)}$ of $G~\wr~S_{(n)}$, the Fourier transform is
  
\[ \widehat{P}(\rho) = \displaystyle \left[ \frac{n_1}{n^2} + \frac{n_1(n_1-1)}{n^2} r(\lambda_1) \right] I, \]

\noindent
where $(n) = (n_1, \ldots, n_s)$, $(\lambda) = ([\lambda_1], \ldots, [\lambda_s])$, and $r(\lambda_1) =
\chi_{[\lambda_1]}(\tau)/d_{[\lambda_1]}$ with transposition $\tau \in S_{n_1}$.
\end{lemma}

\proof{Proof} Recall from Lemma~\ref{3.3.1} that $P$ is constant on conjugacy classes.  It then follows from Lemma \ref{2.5.1} that 
$\widehat{P}(\rho) = C \cdot I$, where $C$ is a constant.  By applying the results from Corollary~\ref{3.3.2} and Lemma~\ref{3.5.1},
we find that

\[ \begin{array}{rcl}

C & = & \displaystyle \frac{1}{|G|n}(1)(1) \ \ + \ \ \sum_{k=2}^s \frac{1}{|G|n^2} \left( n |C_k| \right) \sum_{j=1}^s \left(
\frac{n_j}{n} \right) \frac{\chi_{\rho_j}(g_k)}{d_{\rho_j}} \vspace{1pc} \\

  &   & \displaystyle + \ \ \sum_{k=1}^s \frac{2}{|G|^2n^2} \cdot \frac{n(n-1)}{2} \cdot |G| \cdot |C_k| \sum_{j=1}^s 
\frac{n_j(n_j-1)}{n(n-1)} \cdot \frac{\chi_{\rho_j}(g_k)}{d_{\rho_j}} \cdot r(\lambda_j) \vspace{1pc} \\

  & = & \displaystyle \frac{1}{|G|n} \ \ + \ \ \frac{1}{|G|n^2} \sum_{j=1}^s \frac{n_j}{d_{\rho_j}} \sum_{k=2}^s |C_k| \cdot 
\chi_{\rho_j}(g_k) \vspace{1pc} \\

  &   & \displaystyle + \ \ \frac{1}{|G|n^2} \sum_{j=1}^s \frac{n_j(n_j-1)}{d_{\rho_j}} r(\lambda_j) \sum_{k=1}^s |C_k| \cdot 
\chi_{\rho_j}(g_k).

\end{array} \]

\noindent
When $j = 1$, since $\rho_j$ is the trivial representation of $G$ with $d_{\rho_j} = 1$, we have  $\displaystyle
\sum_{k=1}^s |C_k| \cdot \chi_{\rho_j}(g_k) = \sum_{k=1}^s |C_k| = |G|$ and $\displaystyle \sum_{k=2}^s |C_k| \cdot
\chi_{\rho_j}(g_k) = |G| - |C_1| = |G| - 1$.  When $2 \leq j \leq s$, it follows from Proposition \ref{2.4.4} that $\displaystyle
\sum_{k=1}^s |C_k| \cdot \chi_{\rho_j}(g_k) = 0$ and $\displaystyle \sum_{k=2}^s |C_k| \cdot \chi_{\rho_j}(g_k) = - |C_1| \cdot
\chi_{\rho_j}(g_1) = -d_{\rho_j}$.  Thus

\[ \begin{array}{rcl}

C & = & \displaystyle \frac{1}{|G|n} \ \ + \ \ \frac{n_1(|G|-1)}{|G|n^2} \ \ + \ \ \frac{(n-n_1)(-1)}{|G|n^2} \ \ + \ \
\frac{n_1(n_1-1)r(\lambda_1)|G|}{|G|n^2} \vspace{1pc} \\

  & = & \displaystyle \frac{n_1}{n^2} + \frac{n_1(n_1-1)}{n^2} r(\lambda_1). \ \ \ \qed

\end{array} \]

By applying the results from Lemmas~\ref{3.5.1} and \ref{3.6.1} to Lemma~\ref{2.6.1}, we determine all the eigenvalues of the
transition matrix \textbf{P} induced by the probability measure $P$, together with their multiplicities.

\begin{corollary} \label{3.6.2}
Let $P$ be the probability measure on $G~\wr~S_n$ defined in (\ref{3.1.1}).  Let \emph{\textbf{P}} be the transition matrix of the
Markov chain induced by the probability measure $P$.  Then, for the irreducible representation $\rho$ of $G~\wr~S_n$ induced by the
representation $\rho_{(n)} \otimes \rho_{(\lambda)}$ of $G~\wr~S_{(n)}$, there is an eigenvalue $\pi_{\rho}$ of \emph{\textbf{P}}
occurring with algebraic multiplicity 

\[ \displaystyle {n \choose n_1, \ldots, n_s}^2 d_{\rho_1}^{2n_1} \cdots d_{\rho_s}^{2n_s} \cdot d_{[\lambda_1]}^2 \cdots
d_{[\lambda_s]}^2 \]

\noindent
such that

\[ \pi_{\rho} \ \ = \ \ \displaystyle \frac{n_1}{n^2} + \frac{n_1(n_1-1)}{n^2} r(\lambda_1), \]

\noindent
where $(n) = (n_1, \ldots, n_s)$, $(\lambda) = ([\lambda_1], \ldots, [\lambda_s])$, and $r(\lambda_1) = 
\chi_{[\lambda_1]}(\tau)/d_{[\lambda_1]}$ with transposition $\tau \in S_{n_1}$.
\end{corollary}

We have now established the results necessary to prove Theorem~\ref{3.1.3}.

\proof{Proof of Theorem~\ref{3.1.3}} By applying the results from Lemmas~\ref{3.5.1} and \ref{3.6.1} to the Upper Bound Lemma
(\ref{2.6.2}), we find that

\[ \displaystyle \| P^{*k} - U \|_{\mbox{\rm \scriptsize TV}}^2 \ \ \leq \ \ \mbox{$\frac{1}{4}$} \left(|G|^n n!\right) \| P^{*k} - U
\|_2^2 \ \ = \ \ \mbox{$\frac{1}{4}$} \sum_{\rho} d_{\rho}^2 \left[ \frac{n_1}{n^2} + \frac{n_1(n_1-1)}{n^2} r(\lambda_1)
\right]^{2k}, \]

\noindent
where the sum is taken over all nontrivial irreducible representations of $G~\wr~S_n$.  Thus we have

\[ \begin{array}{l}

\| P^{*k} - U \|_{\mbox{\rm \scriptsize TV}}^2 \ \ \leq \ \ \displaystyle \mbox{$\frac{1}{4}$} \left(|G|^n n!\right) \| P^{*k} - U
\|_2^2 \vspace{1pc} \\

\ \ \ = \ \ \displaystyle \mbox{$\frac{1}{4}$} \sum_{(n)} \sum_{(\lambda)} {n \choose n_1, \ldots, n_s}^2
d_{\rho_1}^{2n_1} \cdots d_{\rho_s}^{2n_s} \cdot d_{[\lambda_1]}^2 \cdots d_{[\lambda_s]}^2 \left[ \frac{n_1}{n^2} +
\frac{n_1(n_1-1)}{n^2} r(\lambda_1) \right]^{2k} \vspace{1pc} \\  

\ \ \ = \ \ \displaystyle \mbox{$\frac{1}{4}$} \sum_{n_1=0}^n {n \choose n_1}^2 d_{\rho_1}^{2n_1} \sum_{[\lambda_1]}
d_{[\lambda_1]}^2 \left[ \frac{n_1}{n^2} + \frac{n_1(n_1-1)}{n^2} r(\lambda_1) \right]^{2k} \vspace{1pc} \\

\ \ \ \ \ \ \ \ \displaystyle \times \ \ \sum_{(n_2, \ldots, n_s)} {n-n_1 \choose n_2, \ldots, n_s}^2 d_{\rho_2}^{2n_2} \cdots
d_{\rho_s}^{2n_s} \sum_{([\lambda_2], \ldots, [\lambda_s])} d_{[\lambda_2]}^2 \cdots d_{[\lambda_s]}^2.

\end{array} \]

Consider the direct product $S_{n_2} \times \cdots \times S_{n_s}$.  It follows from Proposition \ref{2.3.5} and Lemma \ref{2.3.2}
that the sum of the squares of the dimensions of all the irreducible representations of $S_{n_2} \times \cdots \times S_{n_s}$ is

\[ \displaystyle \sum_{([\lambda_2], \ldots, [\lambda_s])} d_{[\lambda_2]}^2 \cdots d_{[\lambda_s]}^2 = |S_{n_2} \times \cdots
\times S_{n_s}| = n_2! \cdots n_s!. \]

\noindent
Using the multinomial theorem and Lemma \ref{2.3.2}, we also have

\[ \begin{array}{rcl}

\displaystyle \sum_{(n_2, \ldots, n_s)} {n-n_1 \choose n_2, \ldots, n_s} d_{\rho_2}^{2n_2} \cdots d_{\rho_s}^{2n_s} & = &
\displaystyle \left( d_{\rho_2}^2 + \cdots + d_{\rho_s}^2 \right)^{n-n_1} \vspace{1pc} \\

& = & \left( |G| - d_{\rho_1}^2 \right)^{n-n_1} = (|G|-1)^{n-n_1}.

\end{array} \]

\noindent
Thus

\[ \displaystyle \sum_{(n_2, \ldots, n_s)} {n-n_1 \choose n_2, \ldots, n_s}^2 d_{\rho_2}^{2n_2} \cdots d_{\rho_s}^{2n_s} 
\sum_{([\lambda_2], \ldots, [\lambda_s])} d_{[\lambda_2]}^2 \cdots d_{[\lambda_s]}^2 = (|G|-1)^{n-n_1} (n-n_1)!, \]

\noindent
which combines with the results above to give
\num \begin{equation} \label{3.6.3}
\begin{array}{l}

\| P^{*k} - U \|_{\mbox{\rm \scriptsize TV}}^2 \ \ \leq \ \ \displaystyle \mbox{$\frac{1}{4}$} \left(|G|^n n!\right) \| P^{*k} - U
\|_2^2 \vspace{1pc} \\

\ \ \ = \ \ \displaystyle \mbox{$\frac{1}{4}$} \sum_{n_1=0}^n {n \choose n_1} \frac{n!}{n_1!} (|G|-1)^{n-n_1}
\sum_{[\lambda_1]} d_{[\lambda_1]}^2 \left[ \frac{n_1}{n^2} + \frac{n_1(n_1-1)}{n^2} r(\lambda_1) \right]^{2k},  

\end{array}
\end{equation}

\noindent
where the inner sum is taken over all partitions $[\lambda_1]$ of $n_1$.  The trivial representation of $G~\wr~S_n$, which should
not be included in the summations above, occurs when $n_1 = n$ and $[\lambda_1] = [n]$ (which gives the trivial representation of 
$S_n$).   

For each $1 \leq n_1 \leq n$, it follows from (\ref{2.7.2}) that we may bound the inner sum above, with $[\lambda_1] = [n_1]$
excluded, by

\[ \begin{array}{rcl}

\displaystyle \sum_{[\lambda_1]} d_{[\lambda_1]}^2 \left[ \frac{n_1}{n^2} + \frac{n_1(n_1-1)}{n^2} r(\lambda_1) \right]^{2k}
& = & \displaystyle \left(\frac{n_1}{n}\right)^{4k} \sum_{[\lambda_1]} d_{[\lambda_1]}^2 \left[ \frac{1}{n_1} +
\frac{n_1-1}{n_1} r(\lambda_1) \right]^{2k} \vspace{1pc} \\  

& \leq & \displaystyle \left(\frac{n_1}{n}\right)^{4k} 4a^2 e^{-4c}  

\end{array} \]

\noindent
for a universal constant $a > 0$, when $k \geq \frac{1}{2} n_1\log n_1 + cn_1$.  Since $n \geq n_1$ and $|G| \geq 2$, this is also
true when $k \geq \frac{1}{2} n\log n + \frac{1}{4} n\log(|G|-1) + cn$.   

We must also bound the term for the trivial representation $[\lambda_1] = [n_1]$ for $1 \leq n_1 \leq n-1$.  Since in these cases,
$d_{[\lambda_1]}^2 = 1$ and $r(\lambda_1) = 1$, we have

\[ \displaystyle d_{(\lambda_1)}^2 \left[ \frac{n_1}{n^2} + \frac{n_1(n_1-1)}{n^2} r(\lambda_1) \right]^{2k} = \left[ 
\frac{n_1}{n^2} + \frac{n_1(n_1-1)}{n^2} \right]^{2k} = \left(\frac{n_1}{n}\right)^{4k}. \]

\noindent
These results lead to the upper bound

\[ \begin{array}{rcl}

\| P^{*k} - U \|_{\mbox{\rm \scriptsize TV}}^2  &  \leq  &  \displaystyle \mbox{$\frac{1}{4}$} \left(|G|^n n!\right) \| P^{*k} - U
\|_2^2 \vspace{1pc} \\

& = & \displaystyle \mbox{$\frac{1}{4}$} \sum_{n_1=0}^n {n \choose n_1} \frac{n!}{n_1!} (|G|-1)^{n-n_1} 
\sum_{[\lambda_1]} d_{[\lambda_1]}^2 \left[ \frac{n_1}{n^2} + \frac{n_1(n_1-1)}{n^2} r(\lambda_1) \right]^{2k} \vspace{1pc} \\  

& \leq & \displaystyle a^2 e^{-4c} \sum_{n_1=1}^n {n \choose n_1} \frac{n!}{n_1!} (|G|-1)^{n-n_1} \left(\frac{n_1}{n}\right)^{4k} 
\vspace{1pc} \\

&   & \displaystyle + \ \ \mbox{$\frac{1}{4}$} \sum_{n_1=1}^{n-1} {n \choose n_1} \frac{n!}{n_1!} (|G|-1)^{n-n_1}
\left(\frac{n_1}{n}\right)^{4k}.

\end{array} \]

\noindent
Now notice that, when $k = \frac{1}{2}n\log n + \frac{1}{4}n\log(|G|-1) + cn$, then

\[ \displaystyle \left(\frac{n_1}{n}\right)^{4k} \ = \ \left(\frac{n_1}{n}\right)^{-n \left[-2\log(n) - \log(|G|-1) - 4c\right]} 
\ = \ \left[\frac{e^{-4c}}{(|G|-1)n^2}\right]^{-n \log\left( n_1/n \right)}, \]

\noindent
which combines with the results above to give

\[ \begin{array}{rcl}

\| P^{*k} - U \|_{\mbox{\rm \scriptsize TV}}^2  &  \leq  &  \displaystyle \mbox{$\frac{1}{4}$} \left(|G|^n n!\right) \| P^{*k} - U
\|_2^2 \vspace{1pc} \\

& \leq & \displaystyle a^2 e^{-4c} \sum_{n_1=1}^n {n \choose n_1} \frac{n!}{n_1!} (|G|-1)^{n-n_1} 
\left[\frac{e^{-4c}}{(|G|-1)n^2}\right]^{-n \log\left( n_1/n \right)} \vspace{1pc} \\  

&   & \displaystyle + \ \ \mbox{$\frac{1}{4}$} \sum_{n_1=1}^{n-1} {n \choose n_1} \frac{n!}{n_1!} (|G|-1)^{n-n_1} 
\left[\frac{e^{-4c}}{(|G|-1)n^2}\right]^{-n \log\left( n_1/n \right)}.  

\end{array} \]

\noindent
If we let $i = n-n_1$, it follows that

\[ \begin{array}{rcl}

\| P^{*k} - U \|_{\mbox{\rm \scriptsize TV}}^2  &  \leq  &  \displaystyle \mbox{$\frac{1}{4}$} \left(|G|^n n!\right) \| P^{*k} - U
\|_2^2 \vspace{1pc} \\

& \leq & \displaystyle a^2 e^{-4c} \sum_{i=0}^{n-1} \frac{1}{i!} \left(e^{-4c}\right)^i \ \ \ + \ \ \ \mbox{$\frac{1}{4}$} e^{-4c}
\sum_{i=0}^{n-2} \frac{1}{(i+1)!} \left(e^{-4c}\right)^i \vspace{1pc} \\

& \leq & \displaystyle a^2 e^{-4c} \exp\left( e^{-4c} \right) \ \ \ + \ \ \ \mbox{$\frac{1}{4}$} e^{-4c} \exp\left( e^{-4c} \right).

\end{array} \]

\noindent
Since $c > 0$, we have $\exp(e^{-4c}) < e$.  Therefore

\[ \| P^{*k} - U \|_{\mbox{\rm \scriptsize TV}}^2 \ \ \leq \ \ \displaystyle \mbox{$\frac{1}{4}$} \left(|G|^n n!\right) \| P^{*k} - U
\|_2^2 \ \ \leq \ \ \left[ \left( a^2 + \mbox{$\frac{1}{4}$} \right) e \right] e^{-4c}, \]

\noindent
from which the desired result follows.  $\qed$

Theorem~\ref{3.1.3} shows that $k = \frac{1}{2} n \log n + \frac{1}{4} n \log (|G|-1) + cn$ steps are sufficient for the
(normalized) $\ell^2$ distance, and hence also the total variation distance, to become small.  A lower bound in the (normalized)
$\ell^2$ metric can also be derived by examining $n^2 (|G|-1) \left( 1 - \frac{1}{n} \right)^{4k}$, which is the dominant
contribution to the summation (\ref{3.6.3}) from the proof of Theorem~\ref{3.1.3}.  This term corresponds to the choice $n_1 = n-1$
with $[\lambda_1] = [n-1]$.  Notice that $k = \frac{1}{2} n \log n + \frac{1}{4} n \log (|G|-1) - cn$ steps are necessary for just
this term to become small.

That $k = \frac{1}{2} n \log n - cn$ steps are \emph{necessary} for the total variation distance to become small follows directly from
Theorem \ref{2.7.3}, as follows.  Recall that Theorem \ref{2.7.3} shows that $k = \frac{1}{2} n \log n - cn$ steps are necessary for
the total variation distance to uniformity to become small for a random walk generated by random transpositions from the symmetric
group $S_n$.  This is exactly the random walk on $G~\wr~S_n$ introduced in Section~\ref{3.1}, if the vector $\vec{x}$ from
$(\vec{x};\pi) \in G~\wr~S_n$ is ignored.  Thus Theorem \ref{2.7.3} provides a lower bound on the distance to uniformity in the total
variation metric, and hence also in the (normalized) $\ell^2$ metric.

In the special case of the hyperoctahedral group $\mathbb{Z}_2~\wr~S_n$, we have $|G| = 2$.  Thus the upper bound in 
Theorem~\ref{3.1.3} matches the lower bound derived from Theorem \ref{2.7.3}.

That $k = \frac{1}{2} n \log n + cn$ steps are also \emph{sufficient} for the total variation distance to become small is a result of
the following.

\begin{theorem} \label{3.6.4}
Let $P$ and $U$ be the probability measures on the complete monomial group $G~\wr~S_n$ defined in (\ref{3.1.1}) and
(\ref{3.1.2}), respectively.  Let $k = \frac{1}{2} n \log n + cn$.  Then there exists a universal constant $\hat{b} > 0$ such that

\[ \| P^{*k} - U \|_{\mbox{\rm \scriptsize TV}} \ \ \leq \ \ \hat{b}e^{-2c} \ \ \ \mathrm{for\ all\ } \mbox{$c > 0$}. \]

\end{theorem}

\proof{Proof} The probability measure $P$ on $G~\wr~S_n$ induces probability measures $Q$ on $G^n$ and $R$ on $S_n$ by defining

\[ Q(\vec{x}) \ := \ \sum_{\pi \in S_n} P(\vec{x};\pi) \ \ \ \mathrm{and} \ \ \ \mbox{$R(\pi) \ := \ \displaystyle \sum_{\vec{x} \in
G^n} P(\vec{x};\pi)$}. \]

\noindent
Notice that $R$ is the probability measure on $S_n$ defined in (\ref{2.1.1}).

Recall that $P$ is the one-step distribution for a random walk \textbf{W}$\ = (W_0, W_1, W_2, \ldots)$ on $G~\wr~S_n$, as described
in Section~\ref{2.6}.  Define stochastic processes \textbf{X}$\ = (X_0, X_1, X_2, \ldots)$, with state space $G^n$, and \textbf{Y}$\ 
= (Y_0, Y_1, Y_2, \ldots)$, with state space $S_n$, by setting  $W_k =: (X_k; Y_k)$ for $k \geq 0$.  It is not hard to see
that \textbf{X} (resp., \textbf{Y}) is a random walk on $G^n$ (resp., on $S_n$) with one-step distribution $Q$ (resp., $R$).  

For each $1 \leq i \leq n$, let $T_i$ be the step index $k$ for the random walk \textbf{W} at which the element $x_i \in \vec{x}$ is
first multiplied by a uniformly chosen random element of $G$, as the result of an element either of the form $(\vec{u};e)$ or of
the form $(\vec{v};\tau)$.  Define $\displaystyle T := \max_{1 \leq i \leq n} T_i$.  Thus $T$ is the step at which the last element
of $\vec{x}$ is randomized, and

\[ \mathbb{P} \left\{ T > k \right\} \ \leq \ \sum_{i=1}^n \mathbb{P} \left\{ T_i > k \right\} \ = \ \sum_{i=1}^n \left( 1 -
\mbox{$\frac{1}{n}$} \right)^{2k} \ \leq \ n e^{-2k/n}. \]

\noindent
Notice that

\[ \begin{array}{rcl}

P^{*k}(\vec{x};\pi)  &  \geq  &  \displaystyle \mathbb{P} \left\{ W_k = (\vec{x};\pi), T \leq k \right\} \vspace{1pc} \\

&  =  &  \displaystyle \frac{1}{|G|^n} \mathbb{P} \left\{ Y_k = \pi, T \leq k \right\} \vspace{1pc} \\

&  =  &  \displaystyle \frac{1}{|G|^n} \left[ R^{*k}(\pi) \ - \ \mathbb{P} \left\{ Y_k = \pi, T > k \right\} \right].

\end{array} \]

\noindent
Thus

\[ \begin{array}{l}

\displaystyle \| P^{*k} - U \|_{\mbox{\rm \scriptsize TV}} \ \ = \ \ \displaystyle \max_{A \subseteq G} \left[ U(A) - P^{*k}(A) \right]
\vspace{1pc} \\

\ \ \ \leq \ \ \displaystyle \max_{A \subseteq G} \sum_{(\vec{x};\pi) \in A} \left[ \frac{1}{|G|^n n!} \ - \ \frac{1}{|G|^n}
R^{*k}(\pi) \ + \ \frac{1}{|G|^n} \mathbb{P} \left\{ Y_k = \pi, T > k \right\} \right] \vspace{1pc} \\

\ \ \ \leq \ \ \displaystyle \max_{A \subseteq G} \sum_{(\vec{x};\pi) \in A} \left[ \frac{1}{|G|^n n!} \ - \ \frac{1}{|G|^n}
R^{*k}(\pi) \right] \ + \sum_{(\vec{x};\pi) \in G~\wr~S_n} \frac{1}{|G|^n} \mathbb{P} \left\{ Y_k = \pi, T > k \right\} \vspace{1pc}
\\

\ \ \ = \ \ \displaystyle \| R^{*k} - U_{S_n} \|_{\mbox{\rm \scriptsize TV}} \ + \ \mathbb{P} \left\{ T > k \right\}.

\end{array} \]

Let $k = \frac{1}{2} n \log n + cn$.  It follows from Theorem \ref{2.1.3} that there exists a universal constant $a > 0$ such that
$\| R^{*k} - U_{S_n} \|_{\mbox{\rm \scriptsize TV}} \ \leq \ a e^{-2c}$ for all $c > 0$, where $U_{S_n}$ is the uniform distribution on
$S_n$ as defined in (\ref{2.1.2}).  Furthermore, it follows from above that $\mathbb{P} \left\{ T > k \right\} \ \leq \ e^{-2c}$.
Therefore

\[ \| P^{*k} - U \|_{\mbox{\rm \scriptsize TV}} \ \ \leq \ \ \left( a + 1 \right) e^{-2c}, \]

\noindent
from which the desired result follows. \qed

Notice that, for the independent shuffles random walk, the rate of convergence to uniformity for the (normalized) $\ell^2$ distance
is slightly slower (due to the addition of the term $\frac{1}{4} n \log(|G|-1)$) than that for the total variation distance. 
However, if $|G|$ is moderate relative to $n$, these rates of convergence are ``essentially'' the same.

The following table summarizes the number of steps (both necessary and sufficient) for the distance (both normalized $\ell^2$ and
total variation) to uniformity to become small for various special cases of the independent shuffles random walk analyzed in this
section.

\newpage

\threespace

\begin{center}

\begin{tabular}{||c|c|c|l|r||} \hline

\multicolumn{5}{||c||}{Random walk on $G~\wr~S_n$} \\

\multicolumn{5}{||c||}{(with independent randomizations)} \\ \hline\hline

$G$  &  metric  &  nec.\ or suff.  &  number of steps  &  proof  \\ \hline\hline

  &  $\ell^2$  &  sufficient  &  $\frac{1}{2} n \log n$  &  Thm.\ \ref{3.1.3}  \\ \cline{3-5}

$\mathbb{Z}_2$   &  &  necessary  &  $\frac{1}{2} n \log n$  &  pf.\ of Thm.\ \ref{3.1.3}  \\ \cline{2-5}

  &  $TV$  &  sufficient  &  $\frac{1}{2} n \log n$  &  Thm.\ \ref{3.1.3}  \\ \cline{3-5}

  &  &  necessary  &  $\frac{1}{2} n \log n$  &  Thm.\ \ref{2.7.3}  \\ \hline

  &  $\ell^2$  &  sufficient  &  $\frac{1}{2} n \log n + \frac{1}{4} n \log(m-1)$  &  Thm.\ \ref{3.1.3} \\ \cline{3-5}

$\mathbb{Z}_m$   &  &  necessary  &  $\frac{1}{2} n \log n + \frac{1}{4} n \log(m-1)$  &  pf.\ of Thm.\ \ref{3.1.3}  \\ \cline{2-5}

  &  $TV$  &  sufficient  &  $\frac{1}{2} n \log n$  &  Thm.\ \ref{3.6.4}  \\ \cline{3-5}

  &  &  necessary  &  $\frac{1}{2} n \log n$  &  Thm.\ \ref{2.7.3}  \\ \hline

  &  $\ell^2$  &  sufficient  &  $\frac{1}{2} n \log n + \frac{1}{4} n \log(|m!|-1)$  &  Thm.\ \ref{3.1.3} \\ \cline{3-5}

$S_m$  &  &  necessary  &  $\frac{1}{2} n \log n + \frac{1}{4} n \log(|m!|-1)$  &  pf.\ of Thm.\ \ref{3.1.3}  \\ \cline{2-5}

  &  $TV$  &  sufficient  &  $\frac{1}{2} n \log n$  &  Thm.\ \ref{3.6.4}  \\ \cline{3-5}

  &  &  necessary  &  $\frac{1}{2} n \log n$  &  Thm.\ \ref{2.7.3}  \\ \hline 

  &  $\ell^2$  &  sufficient  &  $\frac{1}{2} n \log n + \frac{1}{4} n \log(|G|-1)$  &  Thm.\ \ref{3.1.3} \\ \cline{3-5}

$G$  &  &  necessary  &  $\frac{1}{2} n \log n + \frac{1}{4} n \log(|G|-1)$  &  pf.\ of  Thm.\ \ref{3.1.3}  \\ \cline{2-5}

abelian  &  $TV$  &  sufficient  &  $\frac{1}{2} n \log n$  &  Thm.\ \ref{3.6.4}  \\ \cline{3-5}

  &  &  necessary  &  $\frac{1}{2} n \log n$  &  Thm.\ \ref{2.7.3}  \\ \hline

  &  $\ell^2$  &  sufficient  &  $\frac{1}{2} n \log n + \frac{1}{4} n \log(|G|-1)$  &  Thm.\ \ref{3.1.3} \\ \cline{3-5}

$G$  &  &  necessary  &  $\frac{1}{2} n \log n + \frac{1}{4} n \log(|G|-1)$  &  pf.\ of Thm.\ \ref{3.1.3}  \\ \cline{2-5}

nonabelian  &  $TV$  &  sufficient  &  $\frac{1}{2} n \log n$  &  Thm.\ \ref{3.6.4}  \\ \cline{3-5}

  &  &  necessary  &  $\frac{1}{2} n \log n$  &  Thm.\ \ref{2.7.3}  \\ \hline 

\end{tabular}

%\vspace{1pc}

%\textbf{Table 1}

\end{center}

\newpage

\twospace

\subsection{Analysis of the Paired Shuffles Random Walk.} \label{3.7}

We now describe a slight variant of the independent shuffles random walk introduced in Section~\ref{3.1}.  This will provide a second 
benchmark random walk, with known rate of convergence, on $G~\wr~S_n$ for use in the comparison technique.  Schoolfield (1998) analyzed
a random walk for which comparisons to the independent shuffles and paired shuffles random walks, in the special case of the hyperoctahedral 
group $\mathbb{Z}_2~\wr~S_n$, gave different bounds.

Imagine $n$ decks of cards, labeled 1 through $n$, in sequential order, each with its $m$ cards in sequential order.  Independently 
choose two integers $p$ and $q$, uniformly from $\{1, 2, \ldots, n\}$.   

If $p \neq q$, transpose the decks in positions $p$ and $q$.  Then, independently of the choice of $p$ and $q$ and uniformly (i.e.,
with probability $\frac{1}{|G|} = \frac{1}{m!}$ each), permute the deck terminating in position $p$ by a permutation $\pi \in G =
S_m$ and permute the deck terminating in position $q$ by $\pi^{-1} \in G = S_m$.  Notice that the only elements of the form
$(\vec{v};\tau)$ that occur in this combination of operations are from the single conjugacy class $C_1^{(\vec{v};\tau)}$.  The
probability that an element from any of the the other $s-1$ conjugacy classes occurs now vanishes.

If $p = q$ (which occurs with probability $1/n$), leave the decks in their current positions.  Then, again independently and 
uniformly, permute the deck in position $p = q$ by a permutation in $G = S_m$.  The probabilities of the identity and of the 
elements $(\vec{u};e)$ are thus unchanged from the independent shuffles random walk.

Again, as in Section~\ref{3.1}, we will actually examine a random walk on $G~\wr~S_n$ for any group $G$, not just the symmetric
group $S_m$.  In this more general case, the example above is equivalent to beginning with a vector $(e, \ldots, e) \in G^n$,
where each $e \in G$ is the identity element.  Two elements of this vector are then transposed as the decks of cards were above.  The
transposed elements of this vector are then multiplied by elements of $G$ as the individual decks were permuted above.

We refer to the process on $G~\wr~S_n$ described above as the \emph{paired shuffles} random walk, again retaining use of the word 
``shuffles'' even when $G$ is not necessarily $S_m$.  As in Section~\ref{3.1}, the \emph{paired spins} and \emph{paired flips} random
walks are defined in an analogous manner in the special cases of the generalized symmetric group $\mathbb{Z}_m~\wr~S_n$ and the
hyperoctahedral group $\mathbb{Z}_2~\wr~S_n$, respectively.

The paired shuffles random walk may be modeled formally by a probability measure $Q$ on the complete monomial group $G~\wr~S_n$.  We
may thus define the following probability measure on the set of all elements of $G~\wr~S_n$:
\num \begin{equation} \label{3.7.1}
\begin{array}{rcll}

Q(\vec{e};e) & = & \displaystyle \frac{1}{|G|n}, & \vspace{.5pc} \\

Q(\vec{u};e) & = & \displaystyle \frac{1}{|G|n^2} & \mbox{where $\vec{u} \neq \vec{e} \in G^n$}, \vspace{.5pc} \\

Q(\vec{v};\tau) & = & \displaystyle \frac{2}{|G| n^2} & \mbox{where $(\vec{v};\tau) \in C_1^{(\vec{v};\tau)}$}, 
\vspace{.5pc} \\

Q(\vec{x};\pi) & = & 0 & \mbox{otherwise},

\end{array}
\end{equation}

\noindent
where there is only one non-identity entry of $\vec{u} \in G^n$, and where if $\tau \in S_n$ is the transposition $(p\ q)$ then the
only possible non-identity entries of $\vec{v} \in G^n$ are in positions $p$ and $q$ and these entries are mutually inverse elements
of $G$.  In the special case of the hyperoctahedral group $\mathbb{Z}_2~\wr~S_n$, the conjugacy class $C_1^{(\vec{v};\tau)}$ is the
even transpositions.

In order to continue our analysis of the paired shuffles random walk, we must now calculate the Fourier transform of $Q$ at each
irreducible representation of the complete monomial group $G~\wr~S_n$.

\begin{lemma} \label{3.7.2}
Let $Q$ be the probability measure on $G~\wr~S_n$ defined in (\ref{3.7.1}).  For the irreducible representation $\rho$ of
$G~\wr~S_n$ induced by the representation $\rho_{(n)} \otimes \rho_{(\lambda)}$ of $G~\wr~S_{(n)}$, the Fourier transform is
  
\[ \widehat{Q}(\rho) = \displaystyle \left[ \frac{n_1}{n^2} + \sum_{j=1}^s \frac{n_j(n_j-1)}{n^2} r(\lambda_j) \right] I, \]

\noindent
where $(n) = (n_1, \ldots, n_s)$, $(\lambda) = ([\lambda_1], \ldots, [\lambda_s])$, and $r(\lambda_j) =
\chi_{[\lambda_j]}(\tau)/d_{[\lambda_j]}$ with transposition $\tau \in S_{n_j}$.
\end{lemma}

\proof{Proof} Notice that $Q$ is constant on the conjugacy classes of $G~\wr~S_n$.  It then follows from Lemma \ref{2.5.1} that
$\widehat{Q}(\rho) = C \cdot I$, where $C$ is a constant.  By applying the results from Corollary~\ref{3.3.2} and Lemma~\ref{3.5.1},
we find that

\[ \begin{array}{rcl}

C & = & \displaystyle \frac{1}{|G|n}(1)(1) \ \ + \ \ \sum_{k=2}^s \frac{1}{|G|n^2} \left( n |C_k| \right) \sum_{j=1}^s \left(
\frac{n_j}{n} \right) \frac{\chi_{\rho_j}(g_k)}{d_{\rho_j}} \vspace{1pc} \\

  &   & \displaystyle + \ \ \frac{2}{|G| n^2} \cdot \frac{n(n-1)}{2} \cdot |G| \sum_{j=1}^s \frac{n_j(n_j-1)}{n(n-1)} r(\lambda_j)
\vspace{1pc} \\

  & = & \displaystyle \frac{1}{|G|n} \ \ + \ \ \frac{1}{|G|n^2} \sum_{j=1}^s \frac{n_j}{d_{\rho_j}} \sum_{k=2}^s |C_k| \cdot 
\chi_{\rho_j}(g_k) \ \ + \ \ \sum_{j=1}^s \frac{n_j(n_j-1)}{n^2} r(\lambda_j) \vspace{1pc} \\

  & = & \displaystyle \frac{1}{|G|n} \ \ + \ \ \frac{n_1(|G|-1)}{|G|n^2} \ \ + \ \ \frac{(n-n_1)(-1)}{|G|n^2} \ \ + \ \ \sum_{j=1}^s
\frac{n_j(n_j-1)}{n^2} r(\lambda_j) \vspace{1pc} \\

  & = & \displaystyle \frac{n_1}{n^2} + \sum_{j=1}^s \frac{n_j(n_j-1)}{n^2} r(\lambda_j). \ \ \ \qed

\end{array} \]

By applying the results from Lemmas~\ref{3.5.1} and \ref{3.7.2} to Lemma~\ref{2.6.1}, we determine all the eigenvalues of the
transition matrix \textbf{Q} induced by the probability measure $Q$, together with their multiplicities.

\begin{corollary} \label{3.7.3}
Let $Q$ be the probability measure on $G~\wr~S_n$ defined in (\ref{3.7.1}).  Let \emph{\textbf{Q}} be the transition matrix of the
Markov chain induced by the probability measure $Q$.  Then, for the irreducible representation $\rho$ of $G~\wr~S_n$ induced by the
representation $\rho_{(n)} \otimes \rho_{(\lambda)}$ of $G~\wr~S_{(n)}$, there is an eigenvalue $\pi_{\rho}$ of \emph{\textbf{Q}}
occurring with algebraic multiplicity 

\[ \displaystyle {n \choose n_1, \ldots, n_s}^2 d_{\rho_1}^{2n_1} \cdots d_{\rho_s}^{2n_s} \cdot d_{[\lambda_1]}^2 \cdots
d_{[\lambda_s]}^2 \]

\noindent
such that

\[ \pi_{\rho} \ \ = \ \ \displaystyle \frac{n_1}{n^2} + \sum_{j=1}^s \frac{n_j(n_j-1)}{n^2} r(\lambda_j), \]

\noindent
where $(n) = (n_1, \ldots, n_s)$, $(\lambda) = ([\lambda_1], \ldots, [\lambda_s])$, and $r(\lambda_j) = 
\chi_{[\lambda_j]}(\tau)/d_{[\lambda_j]}$ with transposition $\tau \in S_{n_j}$.
\end{corollary}

The following result establishes an upper bound on both the total variation distance and the $\ell^2$ distance between $Q^{*k}$ and
the uniform distribution $U$ on $G~\wr~S_n$.  The total variation upper bound is rather poor when $G$ is nonabelian, as shown by
Theorem~\ref{3.7.6}.  The quality of the $\ell^2$ upper bound will be discussed following the proof of the theorem.

\begin{theorem} \label{3.7.4}
Let $Q$ and $U$ be the probability measures on the complete monomial group $G~\wr~S_n$ defined in (\ref{3.7.1}) and
(\ref{3.1.2}), respectively.  Let

\[ k = \max\left\{ n \log n + \mbox{$\frac{1}{2}$} n \log (|G|-1) + \mbox{$\frac{1}{2}$} n \log (s-1), \ \mbox{$\frac{1}{2}$} n \log
\delta_n \right\} + 2cn, \]

\noindent
where $\displaystyle \delta_n := \sum_{j=2}^s d_{\rho_j}^{2n}$ and $d_{\rho_j}$ is the dimension of the irreducible representation
$\rho_j$ of $G$ for $2 \leq j \leq s$.  Then there exists a universal constant $b>0$ such that

\[ \| Q^{*k} - U \|_{\mbox{\rm \scriptsize TV}} \ \ \leq \ \ \mbox{$\frac{1}{2}$} \left(|G|^n n!\right)^{1/2} \| Q^{*k} - U
\|_2 \ \ \leq \ \ be^{-2c} \ \ \ \mathrm{for\ all\ } \mbox{$c > 0$}. \]

\end{theorem}

\proof{Proof} By applying the results from Lemmas~\ref{3.5.1} and \ref{3.7.2} to the Upper Bound Lemma (\ref{2.6.2}), we find that
\num \begin{equation} \label{3.7.5}
\begin{array}{l}

\displaystyle \| Q^{*k} - U \|_{\mbox{\rm \scriptsize TV}}^2 \ \ \leq \ \ \displaystyle \mbox{$\frac{1}{4}$} \left(|G|^n n!\right) \|
Q^{*k} - U \|_2^2 \ \ = \ \ \mbox{$\frac{1}{4}$} \sum_{\rho} d_{\rho}^2 \left[ \frac{n_1}{n^2} + \sum_{j=1}^s
\frac{n_j(n_j-1)}{n^2} r(\lambda_j) \right]^{2k} \vspace{1pc} \\

\ \ \ = \ \ \displaystyle \mbox{$\frac{1}{4}$} \sum_{(n)} \sum_{(\lambda)} {n \choose n_1, \ldots, n_s}^2 d_{\rho_1}^{2n_1} \cdots
d_{\rho_s}^{2n_s} \cdot d_{[\lambda_1]}^2 \cdots d_{[\lambda_s]}^2 \left[ \frac{n_1}{n^2} + \sum_{j=1}^s \frac{n_j(n_j-1)}{n^2}
r(\lambda_j) \right]^{2k}

\end{array}
\end{equation}

\noindent
where the sums are taken over all nontrivial irreducible representations of $G~\wr~S_n$.  Notice that

\[ \begin{array}{l}

\displaystyle \left[ \frac{n_1}{n^2} + \sum_{j=1}^s \frac{n_j(n_j-1)}{n^2} r(\lambda_j) \right]^{2k} \ \ = \ \ \displaystyle \left\{
\left[ \frac{n_1}{n^2} + \frac{n_1(n_1-1)}{n^2} r(\lambda_1) \right] \ + \ \sum_{j=2}^s \frac{n_j(n_j-1)}{n^2} r(\lambda_j)
\right\}^{2k} \vspace{1pc} \\

\ \ \ = \displaystyle \left\{ \left(\frac{n_1}{n}\right)^2 \left[ \frac{1}{n_1} + \frac{n_1-1}{n_1} r(\lambda_1) \right] \ + \
\sum_{j=2}^s \left(\frac{n_j}{n}\right)^2 \frac{n_j-1}{n_j} r(\lambda_j) \right\}^{2k} \vspace{1pc} \\

\ \ \ \leq \displaystyle \max \left\{ \left(\frac{n_1}{n}\right)^{2k} \left[ \frac{1}{n_1} + \frac{n_1-1}{n_1} r(\lambda_1)
\right]^{2k}, \ \max_{2 \leq j \leq s} \left(\frac{n_j}{n}\right)^{2k} \left[ \frac{n_j-1}{n_j} r(\lambda_j) \right]^{2k} \right\} 
\vspace{1pc} \\

\ \ \ \leq \displaystyle \left(\frac{n_1}{n}\right)^{2k} \left[ \frac{1}{n_1} + \frac{n_1-1}{n_1} r(\lambda_1) \right]^{2k} \ + \
\sum_{j=2}^s \left(\frac{n_j}{n}\right)^{2k} \left[ \frac{n_j-1}{n_j} r(\lambda_j) \right]^{2k},

\end{array} \]

\noindent
where the first inequality is due to the fact that $\displaystyle \left( \sum_{j=1}^s \alpha_j x_j \right)^{2k} \leq \max_{1 \leq j
\leq s} x_j^{2k}$ whenever each $\alpha_j \geq 0$ and $\alpha_1 + \cdots + \alpha_s = 1$.  As noted in the proof of Theorem~5 in
Section~D of Chapter~3 of Diaconis (1988), to every representation $(\lambda_j)$ there corresponds a conjugate representation
$(\lambda'_j)$ such that $r(\lambda_j) = -r(\lambda'_j)$.  So we have

\[ \sum_{[\lambda_j]} \left[ \frac{n_j-1}{n_j} r(\lambda_j) \right]^{2k} \ \ \leq \ \ 2 \sum_{[\lambda_j]: r(\lambda_j) \geq 0}
\left[ \frac{n_j-1}{n_j} r(\lambda_j) \right]^{2k}. \]

\noindent
Thus

\[ \begin{array}{rcl}

\displaystyle \sum_{[\lambda_j]} \left[ \frac{n_j-1}{n_j} r(\lambda_j) \right]^{2k}  &  \leq  &  \displaystyle 2 \sum_{[\lambda_j]:
r(\lambda_j) \geq 0} \left[ \frac{1}{n_j} + \frac{n_j-1}{n_j} r(\lambda_j) \right]^{2k} \vspace{1pc} \\

&  \leq  &  \displaystyle 2 \ \sum_{[\lambda_j]} \left[ \frac{1}{n_j} + \frac{n_j-1}{n_j} r(\lambda_j) \right]^{2k},

\end{array} \]

\noindent
which combines with the previous results to give 

\[ \begin{array}{l}

\displaystyle \| Q^{*k} - U \|_{\mbox{\rm \scriptsize TV}}^2 \ \ \leq \ \ \displaystyle \mbox{$\frac{1}{4}$} \left(|G|^n n!\right) \|
Q^{*k} - U \|_2^2 \vspace{1pc} \\

\ \ \ \leq \ \ \displaystyle \mbox{$\frac{1}{2}$} \sum_{(n)} \sum_{(\lambda)} {n \choose n_1, \ldots, n_s}^2 d_{\rho_1}^{2n_1}
\cdots d_{\rho_s}^{2n_s} \cdot d_{[\lambda_1]}^2 \cdots d_{[\lambda_s]}^2 \sum_{j=1}^s \left(\frac{n_j}{n}\right)^{2k} \left[ 
\frac{1}{n_j} + \frac{n_j-1}{n_j} r(\lambda_j) \right]^{2k} \vspace{1pc} \\

\ \ \ \ \ \ \ \ + \ \ \displaystyle \mbox{$\frac{1}{4}$} \sum_{j=2}^s d_{\rho_j}^{2n} \left(\frac{n-1}{n}\right)^{2k},

\end{array} \]

\noindent
where the sum over $(\lambda) = ([\lambda_1], \ldots, [\lambda_s])$ is taken over all partitions of $(n_1, \ldots, n_s)$, except
that we omit the trivial partitions $[\lambda_j] = [n]$ for $1 \leq j \leq s$.  The final term reintroduces the appropriate terms
for $[\lambda_j] = [n]$ for $2 \leq j \leq s$.

Continuing as in the proof of Theorem~\ref{3.1.3}, this may be simplified to

\[ \begin{array}{l}

\displaystyle \| Q^{*k} - U \|_{\mbox{\rm \scriptsize TV}}^2 \ \ \leq \ \ \displaystyle \mbox{$\frac{1}{4}$} \left(|G|^n n!\right) \|
Q^{*k} - U \|_2^2 \vspace{1pc} \\

\ \ \ \leq \ \ \displaystyle \mbox{$\frac{1}{2}$} \sum_{j=1}^s \sum_{n_j=0}^n {n \choose n_j} \frac{n!}{n_j!}
(|G|-d_{\rho_j}^2)^{n-n_j} \sum_{[\lambda_j] \vdash n_j} d_{[\lambda_j]}^2 \left(\frac{n_j}{n}\right)^{2k} \left[ \frac{1}{n_j} +
\frac{n_j-1}{n_j} r(\lambda_j) \right]^{2k} \vspace{1pc} \\

\ \ \ \ \ \ \ \ + \ \ \displaystyle \mbox{$\frac{1}{4}$} \sum_{j=2}^s d_{\rho_j}^{2n} \left(\frac{n-1}{n}\right)^{2k} \vspace{1pc}
\\

\ \ \ \leq \ \ \displaystyle \mbox{$\frac{1}{2}$} \ s \sum_{n_j=0}^n {n \choose n_j} \frac{n!}{n_j!} (|G|-1)^{n-n_j}
\sum_{[\lambda_j] \vdash n_j} d_{[\lambda_j]}^2 \left(\frac{n_j}{n}\right)^{2k} \left[ \frac{1}{n_j} + \frac{n_j-1}{n_j}
r(\lambda_j) \right]^{2k} \vspace{1pc} \\

\ \ \ \ \ \ \ \ + \ \ \displaystyle \mbox{$\frac{1}{4}$} \left(\frac{n-1}{n}\right)^{2k} \delta_n

\end{array} \]

\noindent
where, for $1 \leq j \leq s$, the sum $\displaystyle \sum_{[\lambda_j] \vdash n_j}$ is taken over all partitions $[\lambda_j]$ of
$n_j$, with $[\lambda_j] = [n]$ excluded when $n_j = n$, and where $\displaystyle \delta_n = \sum_{j=2}^s d_{\rho_j}^{2n}$.

Recall from the proof of Theorem~\ref{3.1.3} that, when $k \geq \frac{1}{2}n \log n + \frac{1}{4}n \log~(|G|~-~1)  + c'n$,
we may bound the inner sum above using

\[ \sum_{[\lambda_1]} d_{[\lambda_1]}^2 \left(\frac{n_1}{n}\right)^{2k} \left[ \frac{1}{n_1} + \frac{n_1-1}{n_1} r(\lambda_1) 
\right]^{2k} \ \ \leq \ \ 4a^2 e^{-4c'} \left(\frac{n_1}{n}\right)^{2k} \ + \ \left(\frac{n_1}{n}\right)^{2k}. \]

\noindent
for $1 \leq n_1 \leq n-1$, and using

\[ \sum_{[\lambda_1]} d_{[\lambda_1]}^2 \left(\frac{n_1}{n}\right)^{2k} \left[ \frac{1}{n_1} + \frac{n_1-1}{n_1} r(\lambda_1)
\right]^{2k} \ \ \leq \ \ 4a^2 e^{-4c'} \]

\noindent
for $n_1 = n$.  So this is also true when $k \geq n \log n + \frac{1}{2} n \log (|G|-1) + \frac{1}{2} n \log (s-1) + 2cn$, in which
case we may choose $c' = \frac{1}{2} \log n + \frac{1}{4} \log (|G|-1) + \frac{1}{2} \log (s-1) + 2c$.  Thus

\[ \displaystyle e^{-4c'} \ = \ \frac{e^{-8c}}{(|G|-1) (s-1)^2 n^2}. \]

Now notice that when $k \geq n \log n + \frac{1}{2} n \log (|G|-1) + \frac{1}{2} n \log (s-1) + 2cn$,

\[ \left(\frac{n_1}{n}\right)^{2k} \leq \left[ \frac{e^{-4c}}{(|G|-1) (s-1) n^2} \right]^{-n \log(n_1/n)}. \]

\noindent
These results lead to the upper bound

\[ \begin{array}{l}

\| Q^{*k} - U \|_{\mbox{\rm \scriptsize TV}}^2 \ \ \leq \ \ \displaystyle \mbox{$\frac{1}{4}$} \left(|G|^n n!\right) \| Q^{*k} - U
\|_2^2 \vspace{1pc} \\

\ \ \ \leq \ \ \displaystyle 2 s a^2 \left[ \frac{e^{-8c}}{(|G|-1) (s-1)^2 n^2} \right] \sum_{n_1=1}^n {n \choose n_1}
\frac{n!}{n_1!} (|G|-1)^{n-n_1} \left[ \frac{e^{-4c}}{(|G|-1) (s-1) n^2} \right]^{-n \log(n_1/n)} \vspace{1pc} \\

\ \ \ \ \ \ \ \ + \ \ \displaystyle \mbox{$\frac{1}{2}$} \ s \sum_{n_1=1}^{n-1} {n \choose n_1} \frac{n!}{n_1!} (|G|-1)^{n-n_1}
\left[ \frac{e^{-4c}}{(|G|-1) (s-1) n^2} \right]^{-n \log(n_1/n)} \vspace{1pc} \\

\ \ \ \ \ \ \ \ + \ \ \displaystyle \mbox{$\frac{1}{4}$} \ \delta_n \left[ \frac{e^{-4c}}{(|G|-1) (s-1) n^2} \right]^{-n \log
\left(1 - \frac{1}{n}\right)}.

\end{array} \]

Continuing as in the proof of Theorem~\ref{3.1.3}, we find that

\[ \begin{array}{rcl}

\| Q^{*k} - U \|_{\mbox{\rm \scriptsize TV}}^2  &  \leq  &  \displaystyle \mbox{$\frac{1}{4}$} \left(|G|^n n!\right) \| Q^{*k} - U
\|_2^2 \vspace{1pc} \\

&  \leq  &  \displaystyle 2 a^2 \mbox{$\left[ \frac{s}{(|G|-1) (s-1)^2 n^2} \right]$} e^{-8c} \exp\left( e^{-4c} \right) \ \ \ + \ \
\ \mbox{$\frac{1}{2}$} \left[ \mbox{$\frac{s}{s-1}$} \right] e^{-4c} \exp\left( e^{-4c} \right) \vspace{1pc} \\

&  &  \displaystyle + \ \ \mbox{$\frac{1}{4}$} \ \delta_n \left[ \mbox{$\frac{e^{-4c}}{(|G|-1) (s-1) n^2}$} \right].

\end{array} \]

\noindent
Since $c > 0$, we have $\exp(e^{-4c}) < e$.  It then follows that

\[ \begin{array}{rcl}

\| Q^{*k} - U \|_{\mbox{\rm \scriptsize TV}}^2  &  \leq  &  \displaystyle \mbox{$\frac{1}{4}$} \left(|G|^n n!\right) \| Q^{*k} - U
\|_2^2 \vspace{1pc} \\

&  \leq  &  \displaystyle \left[ \left( 4a^2 + 1 \right) e \right] e^{-4c} \ \ \ + \ \ \ \left[ \frac{\delta_n}{4(|G|-1) (s-1) n^2}
\right] e^{-4c}.

\end{array} \]

Recall that when $G$ is abelian, $\displaystyle \delta_n = \sum_{j=2}^s d_{\rho_j}^{2n} = |G| - 1$ for all $n$; in this case, the
proof is complete.

When $G$ is nonabelian, $\delta_n$ increases exponentially with $n$.  We must thus reexamine the term $\displaystyle
\mbox{$\frac{1}{4}$} \ \delta_n \left(\frac{n-1}{n}\right)^{2k}$.  Notice that when $k \geq \frac{1}{2} n \log \delta_n + 2cn$,

\[ \mbox{$\frac{1}{4}$} \ \delta_n \left(\frac{n-1}{n}\right)^{2k} \ \ \leq \ \ \mbox{$\frac{1}{4}$} \ \delta_n e^{-2k/n} \ \ \leq \
\ \mbox{$\frac{1}{4}$} \ e^{-4c}. \]

\noindent
Therefore, choosing $k = \max\left\{ n \log n + \frac{1}{2} n \log (|G|-1) + \frac{1}{2} n \log (s-1), \ \frac{1}{2} n \log \delta_n
\right\} + 2cn$ completes the proof. \qed

Theorem~\ref{3.7.4} shows that $k = \max \Big\{ n \log n + \frac{1}{2} n \log (|G|-1) + \frac{1}{2} n \log (s-1), 
\frac{1}{2} n \log \delta_n \Big\} + 2cn$ steps are sufficient for the (normalized) $\ell^2$ distance, and hence the total
variation distance, to become small.  When $n \log n + \frac{1}{2} n \log (|G|-1) + \frac{1}{2} n \log (s-1) \geq \frac{1}{2} n \log
\delta_n$ (which is always the case when $G$ is abelian), a lower bound in the (normalized) $\ell^2$ metric can also be derived by
examining

\[ n^2 \left( 1 - \frac{1}{n} \right)^{4k} \sum_{j=2}^s d_{\rho_j}^2 \ \ = \ \ n^2 (|G|-1) \left( 1 - \frac{1}{n} \right)^{4k}, \]

\noindent
which, in this case, is the dominant contribution to the summation (\ref{3.7.5}) from the proof of Theorem~\ref{3.7.4}.  This comes
from summing (over $2 \leq j \leq s$) the terms corresponding to the choice $n_1 = n-1$ with $[\lambda_1] = [n-1]$ and $n_j = 1$
with $[\lambda_j] = [1]$.  Notice that $k = \frac{1}{2} n \log n + \frac{1}{4} n \log (|G|-1) - cn$ steps are necessary for just
this term to become small.

Notice that, in this case, since $s \leq |G|$, our upper [$n \log n + \frac{1}{2} n \log (|G|-1) + \frac{1}{2} n \log (s-1)$] and
lower [$\frac{1}{2} n \log n + \frac{1}{4} n \log (|G|-1)$] bounds on the number of steps required for the (normalized)
$\ell^2$ distance to become small differ by at most a constant factor, but we have not been able to close this gap.  However, this 
gap will pose no problems in the implementation of the comparison technique, since its results are accurate only up to a constant factor.  
Nonetheless, no such gap exists for total variation distance, as will be shown later.

When $n \log n + \frac{1}{2} n \log (|G|-1) + \frac{1}{2} n \log (s-1) \leq \frac{1}{2} n \log \delta_n$, a matching lower bound in
the (normalized) $\ell^2$ metric can also be derived by examining

\[ \left( \sum_{j=2}^s d_{\rho_j}^{2n} \right) \left( \frac{n-1}{n} \right)^{2k} \ \ = \ \ \delta_n \left( 1 - \frac{1}{n} 
\right)^{2k}, \]

\noindent
which, in this case, is the dominant contribution to the summation (\ref{3.7.5}) from the proof of Theorem~\ref{3.7.4}.  This comes
from summing (over $2 \leq j \leq s$) the terms corresponding to the choice $n_j = n$ with $[\lambda_j] = [n]$.  Notice that $k =
\frac{1}{2} n \log \delta_n - cn$ steps are necessary for just this term to become small.

For fixed nonabelian $G$, notice that $d_{\max}^{2n} \leq \delta_n \leq (s-1) d_{\max}^{2n}$, where $\displaystyle d_{\max} :=
\max_{2 \leq j \leq s} d_{\rho_j}$.  So $2n \log d_{\max} \leq \log \delta_n \leq 2n \log d_{\max} + \log (s-1)$.  Thus, in this
case,  we have thereby determined that very nearly $n^2 \log d_{\max}$ steps are necessary and sufficient to make (normalized)
$\ell^2$ distance small for a fixed nonabelian group $G$.  How large is $d_{\max}$?  Since $\displaystyle \sum_{j=2}^s d_{\rho_j}^2
= |G| - 1$, somewhat crude bounds are

\[ \left( \frac{|G|-1}{s-1} \right)^{1/2} \ \ \leq \ \ d_{\max} \ \ \leq \ \ \left( |G|-1 \right)^{1/2}. \]

\noindent
For $G = S_m$ with $m \geq 3$, for example, these bounds are sufficient to show $\log d_{\max} (m) = \frac{1}{2} m \log m - 
\frac{1}{2} m - O(m^{1/2})$ as $m \to \infty$, since $|G| = m!$ and $s = p(m) \sim \frac{1}{4m\sqrt{3}} \exp\left\{ \pi  
\sqrt{\frac{2m}{3}} \right\}$, where the asymptotic formula for the partition function $p(\cdot)$ is due to Hardy and Ramanujan
(1918) (see, e.g., Hall (1986), Section 4.2).

That $k = \frac{1}{2} n \log n - cn$ steps are \emph{necessary} for total variation distance to become small again follows
directly from Theorem \ref{2.7.3}, exactly as in Section~\ref{3.6}.  That $k = \frac{1}{2} n \log n + cn$ steps are also
\emph{sufficient} (at least in continuous time) for total variation distance to become small is the result of Theorem~\ref{3.7.6},
which follows.

All of the random walks studied thus far have been discrete-time random walks.  We now introduce the continuous-time analogue of
a discrete time random walk.  Changing from discrete to continuous time will be advantageous in the proof of Theorem~\ref{3.7.6}.

Suppose that $P$ is a probability measure defined on a finite group $G$.  The \emph{continuized} chain corresponding to $P$ is the
continuous-time Markov chain on $G$ started at the identity $e$ with transition rates

\[ p \left( g,h \right) \ \ = \ \ P \left( h g^{-1} \right) \]

\noindent
for $g,h \in G$ with $g \neq h$.  We denote the distribution of the chain at time $t$ by $P_t$, which is given by

\[ P_t (g) \ \ := \ \ \sum_{k=0}^{\infty} e^{-t} \ \frac{t^k}{k!} \ P^{*k} (g) \ \ \ \mathrm{for\ } \mbox{$g \in G$}. \]

The following result shows that time $t_n = \frac{1}{2} n \log n + cn$ is sufficient, as $n \to \infty$, for the total variation
distance to become small in the (continuous-time analogue of the) paired shuffles random walk.  The result is established only in
the limit because the proof relies on classical random graph results known only (at least to us) in the limit.

\begin{theorem} \label{3.7.6}
Let $Q$ and $U$ be the probability measures on the complete monomial group $G~\wr~S_n$ defined in (\ref{3.7.1}) and
(\ref{3.1.2}), respectively.  Let $Q_t$ be the distribution at time $t$ of the continuized chain corresponding to $Q$.  Let $t_n 
= \frac{1}{2} n \log n + cn$.  Then there exists a universal constant $\hat{b} > 0$ such that

\[ \limsup_{n \longrightarrow \infty} \| Q_{t_n} - U \|_{\mbox{\rm \scriptsize TV}} \ \ \leq \ \ \hat{b}e^{-c} \ \ \ \mathrm{for\ all\ }
\mbox{$c > 0$}. \]

\end{theorem}

\proof{Proof} The probability measure $Q$ on $G~\wr~S_n$ induces a probability measure $R$ on $S_n$ by defining

\[ R(\pi) \ := \ \sum_{\vec{x} \in G^n} Q(\vec{x};\pi). \]

\noindent
Notice that $R$ is the probability measure on $S_n$ defined in (\ref{2.1.1}).

In order for the paired shuffles continuized chain to achieve randomness, not only must $\pi$ in $(\vec{x};\pi)$ be a random
permutation of $S_n$, but also must the entries of $\vec{x} \in G^n$ be (uniformly) random elements of $G$.  Recall that the values
of the entries in positions $p$ and $q$ of $\vec{x}$ are multiplied (on the left) by mutually inverse elements of $G$ when
$\vec{x}$ is multiplied (on the left) by $(\vec{v};\tau) \in C_1^{(\vec{v};\tau)}$ with $\tau = (p\ q)$.  In order to determine the
amount of time needed to randomize the entries of $\vec{x}$, begin with $n$ labeled vertices.  At each time that an element
$(\vec{v};(p\ q)) \in C_1^{(\vec{v};\tau)}$ is generated by $Q$, consider an edge to be generated between positions $p$ and $q$ in a
graph $\Gamma$.  Let $T$ be the time at which the graph $\Gamma$ becomes connected.  

Let $T^*$ be the first time $t > T$ at which $Q$ generates an element either of the form $(\vec{u};e)$ or $(\vec{e};e)$.  So we may
suppose that at time $T^*$, an element $x_i \in \vec{x}$ is multiplied (on the left) by a uniformly chosen random element $g \in G$.
Thus $x_i$ is now randomized.  Since $T^* > T$, it follows that there exist elements $x_{j_1}, \ldots, x_{j_m} \in \vec{x}$, with
$j_1, \ldots, j_m \neq i$ and $m \geq 1$, such that, for $1 \leq \ell \leq m$, at some time $t_{j_\ell} \leq T$, the entry
$x_{j_\ell}$ was multiplied (on the left) by a uniformly chosen random element $h_\ell \in G$ and $x_i$ was multiplied by
$h_\ell^{-1} \in G$. But since $x_i$ is now randomized, it follows that $x_{j_1}, \ldots, x_{j_m}$ are also.  By then considering
the elements ``paired'' with $x_{j_1}, \ldots, x_{j_m}$, and so forth, this argument continues on to show that at time $T^*$ 
\emph{all} of the elements of $\vec{x}$ are randomized.

Since elements either of the form $(\vec{u};e)$ or $(\vec{e};e)$ are generated by $Q$ at exponential rate $\lambda = 1/n$, we have

\[ \mathbb{P} \left\{ T^* - T \leq s \right\} \ \ = \ \ 1 \ - \ e^{-s/n}. \]

\noindent
It then follows from the independence of $T$ and $T^* - T$ that

\[ \begin{array}{rcl}

\displaystyle \mathbb{P} \left\{ T^* > t \right\}  &  =  &  \displaystyle \int_{s=0}^t \mbox{$\frac{1}{n}$} e^{-s/n} \ \mathbb{P}
\left\{ T > t-s \right\} \, ds \vspace{1pc} \\

&  =  &  \displaystyle \int_{u=0}^{\infty} e^{-u} \ \mathbb{P} \left\{ T > t-un \right\} \mathbb{I} \left\{ u \leq t/n \right\} \, du.

\end{array} \]

\noindent

\noindent
Since each element of the form $(\vec{x};(p\ q)) \in C_1^{(\vec{x};\tau)}$, which transposes a particular pair of entries $\{p,
q\}$, is generated by $Q$ at exponential rate $\lambda = 2/n^2$, the indicator (call it $I_{\{p,q\}}(t)$) of the presence of any
given edge $\{p, q\}$ in $\Gamma$ at time $t$ has expectation

\[ 1 \ - \ e^{-2t/n^2}. \]

\noindent
Moreover (and this is the advantage of working in continuous time), the stochastic processes $I_{\{p,q\}}(\cdot)$ are mutually
independent.  Let $t \equiv t_n = \frac{1}{2} n \log n + cn$.  Then

\[ 1 - e^{-2t/n^2} \ \ = \ \ 1 - \exp\left\{-n^{-1} (\log n + 2c) \right\} \ \ \sim \ \ n^{-1} (\log n + 2c) \ \ \ \mathrm{as\ }
\mbox{$n \to \infty$} \]

\noindent
for fixed $c \in \mathbb{R}$.  It then follows from a classical random graph result of Erd\"{o}s and R\'{e}nyi (1959) (see, e.g.,
Graham et al (1995), Chapter 6, Section 5), that

\[ \mathbb{P} \left\{ T > t-un \right\} \ \ = \ \ \mathbb{P} \left\{ T > \mbox{$\frac{1}{2}$} n \log n + (c-u)n \right\} \ \
\longrightarrow \ \ 1 \ - \ \exp\left\{ -e^{-2(c-u)} \right\} \ \ \ \mathrm{as\ } \mbox{$n \to \infty$} \]

\noindent
for fixed $c,u \in \mathbb{R}$.  Thus, by the bounded convergence theorem,

\[ \begin{array}{rcl}

\displaystyle \lim_{n \longrightarrow \infty} \ \mathbb{P} \left\{ T^* > t \right\}  &  =  &  \displaystyle \int_{u=0}^{\infty}
e^{-u} \ \left[ 1 \ - \ \exp\left\{ -e^{-2(c-u)} \right\} \right] \, du \vspace{1pc} \\

&  \leq  &  \displaystyle \int_{u=0}^{c} e^{-u} e^{-2(c-u)} \, du \ \ + \ \ \int_{u=c}^{\infty} e^{-u} \, du \vspace{1pc} \\

&  =  &  e^{-c} \ - \ e^{-2c} \ + \ e^{-c} \ \ \leq \ \ 2 e^{-c}

\end{array} \]

\noindent
for fixed $c \in \mathbb{R}$.  

It follows exactly as in the proof of Theorem~\ref{3.6.3} that

\[ \| Q_t - U \|_{\mbox{\rm \scriptsize TV}} \ \ \leq \ \ \| R_t - U_{S_n} \|_{\mbox{\rm \scriptsize TV}} \ + \ \mathbb{P} 
\left\{ T^* > t \right\} \]

\noindent
for every $t \geq 0$.  With $t_n = \frac{1}{2} n \log n + cn$, a continuous time analogue of Theorem \ref{2.1.3} (which we
have confirmed) asserts that there exists a universal constant $a' > 0$ such that

\[ \| R_t - U_{S_n} \|_{\mbox{\rm \scriptsize TV}} \ \ \leq \ \ a' e^{-2c} \ \ \ \mathrm{for\ all\ } \mbox{$n \geq 1$} 
\mathrm{\ and\ all\ } \mbox{$c > 0$}, \]

\noindent
where $U_{S_n}$ is the uniform distribution on $S_n$ defined by (\ref{2.1.2}).  Therefore,

\[ \limsup_{n \longrightarrow \infty} \ \| Q_{t_n} - U \|_{\mbox{\rm \scriptsize TV}} \ \ \leq \ \ \left( a' + 2 \right) e^{-c}, \]

\noindent
from which the desired result follows. \qed

The following table summarizes the number of steps (both necessary and sufficient) for the distance (both normalized $\ell^2$ and
total variation) to uniformity to become small for various special cases of the paired shuffles random walk analyzed in this section.

\newpage

\threespace

\begin{center}

\begin{tabular}{||c|c|c|l|r||} \hline

\multicolumn{5}{||c||}{Random walk on $G~\wr~S_n$} \\

\multicolumn{5}{||c||}{(with paired randomizations)} \\ \hline\hline

$G$  &  metric  &  nec.\ or suff.  &  number of steps  &  proof  \\ \hline\hline

  &  $\ell^2$  &  sufficient  &  $n \log n$  &  Thm.\ \ref{3.7.4}  \\ \cline{3-5}

$\mathbb{Z}_2$   &  &  necessary  &  $\frac{1}{2} n \log n$  &  pf.\ of Thm.\ \ref{3.7.4}  \\ \cline{2-5}

  &  $TV$  &  sufficient  &  $\frac{1}{2} n \log n \ \ (n \rightarrow \infty)$  &  Thm.\ \ref{3.7.6}  \\ \cline{3-5}

  &  &  necessary  &  $\frac{1}{2} n \log n$  &  Thm.\ \ref{2.7.3}  \\ \hline

  &  $\ell^2$  &  sufficient  &  $n \log n + n \log(m-1)$  &  Thm.\ \ref{3.7.4} \\ \cline{3-5}

$\mathbb{Z}_m$   &  &  necessary  &  $\frac{1}{2} n \log n + \frac{1}{4} n \log(m-1)$  &  pf.\ of Thm.\ \ref{3.7.4}  \\ \cline{2-5}

  &  $TV$  &  sufficient  &  $\frac{1}{2} n \log n \ \ (n \rightarrow \infty)$  &  Thm.\ \ref{3.7.6}  \\ \cline{3-5}

  &  &  necessary  &  $\frac{1}{2} n \log n$  &  Thm.\ \ref{2.7.3}  \\ \hline

  &  &  &  $\max \bigg\{ \frac{1}{2} n \log \delta_n, \ n \log n$  &  \\

  &  &  sufficient  &  \ \ \ \ \ \ \ \ $+ \frac{1}{2} n \log(|m!|-1)$  &  Thm.\ \ref{3.7.4}  \\

$S_m$  &  $\ell^2$  &  &  \ \ \ \ \ \ \ \ $+ \frac{1}{2} n \log(p(m)-1) \bigg\}$  &  \\ \cline{3-5}

  &  &  necessary  &  $\max \bigg\{ \frac{1}{2} n \log \delta_n, \ \frac{1}{2} n \log n$  &  \\ 

  &  &  &  \ \ \ \ \ \ \ \ $+ \frac{1}{4} n \log(|m!|-1) \bigg\}$  &  pf.\ of Thm.\ \ref{3.7.4}  \\ \cline{2-5}

  &  $TV$  &  sufficient  &  $\frac{1}{2} n \log n \ \ (n \rightarrow \infty)$  &  Thm.\ \ref{3.7.6}  \\ \cline{3-5}

  &  &  necessary  &  $\frac{1}{2} n \log n$  &  Thm.\ \ref{2.7.3}  \\ \hline

  &  $\ell^2$  &  sufficient  &  $n \log n + n \log(|G|-1)$  &  Thm.\ \ref{3.7.4} \\ \cline{3-5}

$G$  &  &  necessary  &  $\frac{1}{2} n \log n + \frac{1}{4} n \log(|G|-1)$  &  pf.\ of  Thm.\ \ref{3.7.4}  \\ \cline{2-5}

abelian  &  $TV$  &  sufficient  &  $\frac{1}{2} n \log n \ \ (n \rightarrow \infty)$  &  Thm.\ \ref{3.7.6}  \\ \cline{3-5}

  &  &  necessary  &  $\frac{1}{2} n \log n$  &  Thm.\ \ref{2.7.3}  \\ \hline

  &  &  &  $\max \bigg\{ \frac{1}{2} n \log \delta_n, \ n \log n$  &  \\

  &  &  sufficient  &  \ \ \ \ \ \ \ \ $+ \frac{1}{2} n \log(|G|-1)$  &  Thm.\ \ref{3.7.4}  \\

$G$  &  $\ell^2$  &  &  \ \ \ \ \ \ \ \ $+ \frac{1}{2} n \log(s-1) \bigg\}$  &  \\ \cline{3-5}

nonabelian  &  &  necessary  &  $\max \bigg\{ \frac{1}{2} n \log \delta_n, \ \frac{1}{2} n \log n$  &  \\ 

  &  &  &  \ \ \ \ \ \ \ \ $+ \frac{1}{4} n \log(|G|-1) \bigg\}$  &  pf.\ of Thm.\ \ref{3.7.4}  \\ \cline{2-5}

  &  $TV$  &  sufficient  &  $\frac{1}{2} n \log n \ \ (n \rightarrow \infty)$  &  Thm.\ \ref{3.7.6}  \\ \cline{3-5}

  &  &  necessary  &  $\frac{1}{2} n \log n$  &  Thm.\ \ref{2.7.3}  \\ \hline

\end{tabular}

%\vspace{1pc}

%\textbf{Table 2}

\end{center}

\newpage

\twospace

\section*{Acknowledgments.} 
This paper formed a portion of the author's Ph.D. dissertation in the Department of Mathematical Sciences at the Johns Hopkins
University.  The author wishes to thank his advisor Jim Fill, whose assistance was invaluable, particularly in the proof of
Theorem~\ref{3.7.6}.  The author also wishes to thank Persi Diaconis for initially suggesting that he analyze a random walk on
the hyperoctahedral group.  It is from that initial challenge that this paper has evolved.

\twospace

\Line{\hfill
\AOPaddress{Clyde H. Schoolfield, Jr.\\
Department of Statistics\\
Harvard University\\
One Oxford Street\\
Cambridge, Massachusetts 02138\\
e-mail:  {\tt clyde@stat.harvard.edu}}}

\end{document}